\documentclass[aop]{imsart}
\usepackage{amsfonts} 

\RequirePackage{amsthm,amsmath}  % This is required for by \numberwithin{equation}{section} and

\usepackage{mathrsfs}
\usepackage[capitalize]{cleveref}

\usepackage{tikz}
\usepackage{amssymb}

\long\def\/*#1*/{}  %% Usage:  \/******************************************* commented out: ...  ********************/

\startlocaldefs
\numberwithin{equation}{section}
\theoremstyle{plain}

\newtheorem{theorem}{Theorem}[section]
\newtheorem{lemma}[theorem]{Lemma}
\newtheorem{corollary}[theorem]{Corollary}
\theoremstyle{definition}

\newtheorem{proposition}{Proposition}[section]
\endlocaldefs

\begin{document}

\begin{frontmatter}
\title{A Simpler Proof of Kakutani's Conjecture on Random Subdivision and Its Generalizations %\thanksref{}
      }
\runtitle{%Asymptotic Behavior of Spacings Splitting Interval with any pdf in ...
           Simple proof of Kakutani's conjecture on random subdivision
         }
%\thankstext{T1}{Footnote to the title with the ``thankstext'' command.}

\begin{aug}
\author{\fnms{Changqing} \snm{Liu}\ead[label=e1]{C.Liu@ClinTFL.com}}  
\address{ClinTFL Ltd. \\
         \printead{e1}\\
}
\end{aug}

\begin{abstract} 
 
Shifting the perspective from division points to division spacings leads not only to a substantially simpler proof of the conjecture but also to stronger results: limiting distributions are obtained for both partition points and spacings, with large-deviation error bounds, and the limiting empirical c.d.f.s are shown to be independent of the splitting scheme.

The approach is naturally extended to arbitrary distributions of division points, with the limiting spacing distribution established. It is shown stationarity of the spacing distribution iff stationarity of the mean spacing iff uniform distribution of the partition points. This means that the results of Dean and Majumdar (2002) and Janson and Neininger (2008) concerning stationarity of mean spacing apply directly to that of spacing distribution. 

The framework provides subtree-size profiles for arbitrary random trees. A new phenomenon is revealed: the profile may undergo transitions at a series of critical points on different asymptotic scales in \(n\).

\end{abstract}


\end{frontmatter}

\section{Introduction}
Fragmentation is a stochastic process, which starts with unit interval and splits it into certain number of subintervals. Fragmentation can be basically classified into two types, the parallel and the sequential. Parallel Fragmentation is a stochastic process where an object splits into multiple parts simultaneously. Let $F_n$ denote the empirical distribution function for a sample of $n$ random variables $X_1, X_2, ..., X_n$ in the interval $(0, 1)$, i.e.

$$
                                   F_n (x) = \frac{1}{n} \sum_{i=1}^n  \mathbf{1}_{\{X_i \le x\}}  
$$
where $\mathbf{1}_{\{X_i \le x\}}$ is the indicator function. In parallel fragmentation it is extensively studied when $X_1, X_2, ..., X_n$ is i.i.d distribution. The equidistribution of the partition points is proved, i.e. 
\begin{equation}\label{Fn=x}
                          F_n (x) = x \, + \,  \overset{s}{\sim}  \frac{1}{\sqrt{n}}   
\end{equation}  
Here  ``$\overset{s}{\sim}$''  denotes a sub-Gaussian error at the indicated scale. Specifically, 
\[
             \mathbb{P}\left( \big| F_{n}(x) - x \big| \ge \frac{\lambda}{\sqrt{n}} \right) \leq C e^{-c {\lambda^2}} 
\]   
where $C$ and $c$ are two unspecified constants \cite{DKW}. Consequently, for every $A > 0$ there is a constant $C_A > 0$ such that 
\[ 
  %              _______
  %             /  A + 1
  %      C_A = / --------
  %           ?     c
       \mathbb{P}\left (|F_{n}(x) - x | \geq C_A \sqrt{\frac{\log n} n} \right) \le n^{-A}       
\]  
In particular,  $F_n(x) - x = O_{\mathbb{P}}\big(n^{-1/2}\big)$.

Kakutani's interval splitting scheme (hereafter, K-model or Kakutani model) is a sequential process. In 1973, physicist Araki (Huzihiro Araki) asked S. Kakutani (the author of Kakutani's fixed point theorem) a question: whether random splitting of the largest subinterval at each stage almost always yielded a uniform distribution of points. This is later referred as Kakutani's conjecture on random subdivision of largest intervals. The conjecture circulated at a joint Statistics and Probability meeting in 1976 at Oberwolfach\footnote{An invitation-only mathematics research workshop, about 50 experts in probability theory or analysis.}  and challenged the audience for a proof \cite{Beran&Fisher2009}. In this pseudo-Olympiad of mathematics, the ``medallists'' are van Zwet \cite{Zwet1978}, Slud \cite{Slud1978, Slud1982} and Lootgieter \cite{Lootgieter1977, Lootgieter1978}, so to speak. Amongst them, van Zwet's proof is hailed as simple and elegant, and was his personal favourite \cite{Beran&Fisher2009, Mason2012}. Employing van Zwet's (1978) method, Pyke (1980) and Pyke \& van Zwet (2004) studied the empirical distribution of spacings. They showed that the empirical distribution of the normalized spacings converges to Uniform$ (0, 2)$, with the deviations converging in distribution to a Brownian bridge \cite{Pyke1980, Pyke&Zwet2004}. 

In this paper we draw attention to a simpler approach that seems to have been in their blind spot. Its simplicity lies not only in the brevity of the proof but also in its accessibility to students of science and engineering without specialised mathematical training. Besides, the proof covers broader results: it establishes results for empirical distributions of both partition points and spacings. In addition, the result is stronger in the sense that the error terms are controlled by large deviation bounds. Most importantly, the method is general and applies to arbitrary distributions of split points and splitting schemes. We obtain the limiting distributions of spacings for arbitrary distribution of split points, including the general case of discrete distributions, while the literature has only addressed Kakutani model -- binary splitting by a uniformly random point.
 
% RECENT DEVELOPMENTS ? IN STATISTICS, Proceedings Grenoble 6?11 September 1976, edited by J.?R. Barra et al. This was a large, scheduled international statistics conference with hundreds of participants and formally published proceedings (Recent Developments in Statistics, North-Holland, 1977).

% Location: Oberwolfach, Germany (Mathematisches Forschungsinstitut Oberwolfach, MFO) 
% Dates: March or October 1976 (exact meeting for this conjecture still under archival search)
% Organized by: The MFO and a small group of invited organizers
% Nature: A small, invitation-only mathematics research workshop, likely in probability theory or analysis
% No formal proceedings for this particular conjecture presentation ? it was informally shared by Richard Dudley at the meeting.
 
%Acknowledgement. The author recalls with pleasure the 1976 stochastics meeting at Oberwolfach where  R. M. Dudley introduced the participants to Kakutani's conjecture and proceeded to shoot down our combined attempts at solving the problem. 

%Kakutani, S.: A problem of equidistribution on the unit interval [0,1]. Proceedings of the Oberwolfach conference on Measure Theory (1975). Lecture Notes in Math. 541. Springer 1976.  American Mathematical Society
 
%  M. Pollicott* and B. Sewell, University of Warwick, March 17, 2021  "An infinite interval version of the ?-Kakutani equidistribution problem"

\section {Notation}
Let $(V_{1}, V_{2}, ... V_{b})$ be the $b$ spacings induced by $(b-1)$ random division points in $(0, 1)$. Define a random function
\begin{equation}\label{E:xi}
        \xi(x) = \sum_{i=1}^{b}  \mathbf{1}_{\{V_{i}<x\}}
\end{equation}
This is a random variable that gives the count of spacings shorter than $x$. Let $F(x)$ denote the expected count function of $\xi$ on $(0, 1)$, i.e.
\begin{equation}\label{E:Hx}
           \mathbb{E} (\xi(x)) = F(x)
%    C_{s}  =  K_{M} \frac{\mu/\mu_{x}}{1-\mu/\mu_{x}} \label{cs}
\end{equation}
In the Kakutani model, the largest subinterval splits into two parts each time, by a division point uniformly distributed in that subinterval, where $b =2$, $F(x) = 2x$. In the case of $b$ uniformly distributed division points, $F(x) = b\big(1-(1-x)^{b-1}\big)$.

In the process under consideration, at time $i$, $i=1, 2, ..., n$, the largest subinterval---whose length is $M_i$---is split. For simplicity, we use $M_i$ for both the subinterval and its size without confusion. By time $n$, the cumulative number of spacings shorter than $ x $ is given by  
\[
              \xi_{1}(x) + \xi_{2}(x) + \hdots + \xi_{n}(x) 
\]
with $\xi_{1}(x)$ being the random variable defined in~(\ref{E:xi}), resulting from splitting of $M_0 = (0,1)$; $\xi_{2}(x)$ resulting from splitting of $M_{1}$,
\begin{equation}
        \xi_2(x) = \sum_{i=1}^{b}  \mathbf{1}_{  \{ V_{i} < x \}  } 
\end{equation}
where $V_i$ are subintervals induced by $(b-1)$ random division points in $M_{1}$; and so on; i.e. $\xi_i(x)$ counts spacings $ < x$  generated at the $i$-th splitting step, analogous to (\ref{E:xi}) but applied to $M_{i-1}$. At each time step, the division points are drawn from the same distribution. Letting $x_c > 0$ be cut-off such that subintervals smaller than it will not be split further. Thus, the process stops if all the subintervals $M_i \leq x_c$. By~(\ref{E:Hx}), we have, for $x \leq x_c$,
\begin{equation}\label{E:Hxi}
              \mathbb{E}_{i-1}(\xi_{i}(x)) = F(x/M_{i-1})
\end{equation} where $E_{i-1}$ is conditional expectation given first $i-1$ subinterval splittings, and $x / M_{i-1} < 1$ because $x \leq x_c < M_{n-1} \leq ... \leq M_{1} \leq M_{0}=1$.

Now write 
\[  
      X_{i} := \xi_{i}(x) - \mathbb{E}_{i-1}(\xi_{i}(x))   
\]
It is easy to check that, given the first $i-1$ subinterval splittings, the conditional expectation of
$X_{i}$ is zero. In notation, 
\begin{equation}
               \mathbb{E}(X_{i} \mid M_{i-1}, ..., M_{1}, M_{0}) = 0   
\end{equation}             
Therefore, 
                  $$S_{n} = X_{1} + X_{2} + ... + X_{n}$$
is a martingale,  $E(S_{n}) = 0 $ and $X_{i}$ is bounded by a constant $c$\,; e.g., in Kakutani model, $|X_{i}| \leq 2$. Azuma--Hoeffding inequality then gives 
\begin{equation} \label{E:Azuma}
       \mathbb{P}(|S_{n} - 0 | \geq \lambda ) \leq 2e^{ -\frac{ \lambda^2}{2nc^2}
                                        } 
\end{equation}
We write~(\ref{E:Azuma}) as
\begin {equation}\label{E:empirical}
      \sum_{i=1}^{n} \xi_{i}(x) = \sum_{i=1}^{n} \mathbb{E}_{i-1}(\xi_{i}(x)) \,+\, \,\overset{s}{\sim}\sqrt{n}  \hspace{1em} 
\end{equation} 
where $ \mathbb{E}_{i-1}(\xi_{i}(x))$ is given by~(\ref{E:Hxi}). That is to say, with high probability (w.h.p.), the difference between
\[
           \xi_{1}(x) + \xi_{2}(x) + \cdots + \xi_{n}(x)
\] 
and 
\[ 
        \mathbb{E}\big(\xi_{1}(x)\big) + \mathbb{E}_{1}\big(\xi_{2}(x)\big) + \cdots + \mathbb{E}_{n-1}\big( \xi_{n}(x) \big) 
\]
is $ \, \overset{s}{\sim}\sqrt n$, so that the deviations much beyond $\sqrt{n \log n}$ have arbitrarily small probability; alternatively, one may state that deviations of order $n^p$ for any $p > 1/2$ have exponentially small probability.

For binary interval splitting process where $b=2$, the total number of spacings is $n+1$ and the empirical distribution function of spacings is
\[
        \frac{1}{n+1} \sum_{i=1}^{n} \xi_{i}(x)
\]
for $x < x_c$. In general case, the empirical c.d.f. is\;
$ \displaystyle
          \frac{1}{1 + n(b-1)} \sum_{i=1}^{n} \xi_{i}(x).
$

\section {Proof of Kakutani's conjecture}
\begin{theorem}\label {E:th_spacing} 
  In the model of Kakutani interval splitting, the empirical c.d.f. of spacings satisfies
\[
          \frac{x}{x_c} \, + \, \overset{s}{\sim}\frac {1}{\sqrt n}  
\]
for $x \leq x_c$.

\end{theorem}

\begin{proof} From (\ref{E:empirical}), the empirical c.d.f. spacings 
\begin{equation}\label{E: Kaku PROOF}
                \frac{1}{n+1} \sum_{i=1}^{n} \xi_{i}(x) = \frac{1}{n+1} \sum_{i=1}^{n} F\left(\frac{x}{M_{i-1}}\right) 
                 \, + \, \overset{s}{\sim} 
                \frac {1}{\sqrt n} \hspace{1em}  
\end{equation} 
The meaning of ``$ \displaystyle \;  \overset{s}{\sim}  \frac {1} {\sqrt n} $'' is defined in (\ref{Fn=x}).  
In Kakutani model, $\displaystyle F \left( \frac{x}{M_{i-1}} \right) = \frac{2x}{M_{i-1}}$, it follows that 
\begin{equation}\label{E: sumM}
       \frac{1}{n+1} \sum_{i=1}^{n} \xi_{i}(x) = \frac{1}{n+1} \sum_{i=1}^{n} \frac{2x}{M_{i-1}} \, + \, 
       \overset{s}{\sim} \frac {1}{\sqrt n}
\end{equation}
When $x=x_c$, all of the $n+1$ subintervals, and only these, are smaller than $x_c$, each contributing 1 to the total count $\sum \mathbf{1}_{\{V_{i}<x_c\}}$. Therefore the left-hand side of (\ref{E: sumM}) equals 1, yielding 
\[
    \frac{1}{n+1} \sum_{i=1}^{n} \frac {2}{M_{i-1}} = \frac {1}{x_c}  (1  \, + \, \overset{s}{\sim} \frac {1}{\sqrt n})  
\] 
Replacing this in~(\ref{E: sumM}), we have   
\begin{align} \nonumber
       \frac{1}{n+1} \sum_{i=1}^{n} \xi_{i}(x) &= \frac {x}{x_c}  (1  \, + \, \overset{s}{\sim} \frac {1}{\sqrt n})   
                                                              \; + \, \overset{s}{\sim} \frac {1}{\sqrt n}
       \\ \nonumber
                                               &= \frac {x}{x_c}  \, + \, \overset{s}{\sim} \frac {1}{\sqrt n}
\end{align}
the desired result. 
\end{proof}
\vspace{-0.5em}

\textbf{Remark.}
\begin{itemize}
  \item   Thus, the mean spacing is $ \displaystyle \frac{x_c}{2}$\,. Therefore, the distribution of the normalized spacings is $\mathrm{Uniform}(0, 2)$, and the total number of empirical partition points converges to $ \displaystyle \frac {2}{x_c} $\,.
\end{itemize}

\textbf{Empirical partition points.}  To study the distribution of partition points, we consider spacings within the sub-interval $(0, y) \subset (0,1)$. $M_{i-1}$ is called the parent interval of $\xi_i(x)$ (and also of $X_i$).  We will show that the empirical spacings in $(0,y)$ converge the same c.d.f. as in Theorem \ref{E:th_spacing}; thus Theorem \ref{E:th_spacing} holds for any interval contained in $(0,1)$, and therefore the empirical distribution of partition points is Uniform(0,1). 

Let $\{ M_{l_i-1} \}_{i=1}^m$ be all the parent subintervals that either contain the point $y$ (called $y$-subintervals) or occur within $(0,y)$ during the fragmentation process. For the $y$-interval $(0, y)$, define
\begin{equation}
         \xi(x) = \sum_{i=1}^{b}  \mathbf{1}_{  \{ V_{i} < x \  \mbox{and} \ V_{i} \subseteq (0, y) \}  }  
        \hspace{1em} 
\end{equation}
i.e., this $\xi(x)$ only counts those $V_i$ that fall inside $(0, y)$.\footnote{A different symbol -- say, $\eta$ -- should be used intead of $\xi$. For the sake of simplicity and readability, we adhere to the current $\xi$ notation.}  Again, we use $V_i$ to denote either the subinterval itself or its length depending on the context. In Kakutani model, $E \, \xi(x)  = 2\min(x,y)$. For $y$-subinterval $M_{i-1} = (a_{i-1}, b_{i-1})$, 
\begin{equation*}
        \xi_{i}(x) = \sum_{i=1}^{b}  \mathbf{1}_{  \{ V_{i} < x \  \mbox{and} \ V_{i} \subseteq (a_{i-1}, y) \}  }  
\end{equation*} 
with $V_i$ is subinterval of $M_{i-1}$ induced by the $(b-1)$ division points in $M_{i-1}$. 
  
Let  
\[
       S_m =  X_{l_1}  + \hdots + X_{l_{m-1}} + X_{l_m}  
\]
$S_m$ satisfies the Azuma-Hoeffding concentration inequality (see Appendix~\ref{Azuma-Heoffding inequality}), i.e. 
\begin{align*} 
      \mathbb{P}(|S_m| \ge \lambda) \le  2\exp \Big (- \frac{\lambda^2}{2mc^2}  \Big)
\end{align*}
The following lemma is proved (see Appendix~\ref{k=log N}). 
\begin{lemma}\label{lem:log n} During the splitting process, $K=O(\log n)$ cuts suffice for the $y$-interval to reach the cutoff $x_c$ w.o.p.\footnote{w.o.p. (with overwhelming probability) means with probability at least $1 - n^{-A}$ for every fixed $A$
}, and 
\[
     \frac  K n = \,  \, \overset{s}{\sim} n^{-1/2} 
\]     
\end{lemma}

Similar to \eqref{E: Kaku PROOF}, 
\begin{equation}   
                \frac{1}{m+1} \sum_{i=1}^{m} \xi_{l_i}(x) = \frac{1}{m+1} \sum_{i=1}^{m} 
                \mathbb{E}_{l_i-1} \big( \xi_{l_i}(x) \big)  \, + \, \overset{s}{\sim}\; \frac {1}{\sqrt m} \hspace{1em}  
\end{equation}
By Lemma \ref{lem:log n}, it follows 
\begin{align*} 
        \frac{1}{m+1} \sum_{i=1}^{m} \xi_{l_i}(x) 
        &= \frac{1}{m+1} \sum_{i=1}^{m}  \frac{2x}{M_{l_i-1} } \,  - \,  \frac{O(K)\,
\begingroup \renewcommand{\thefootnote}{\fnsymbol{footnote}}\footnotemark[1] \endgroup 
        }{m+1}
          \, + \, \overset{s}{\sim} \frac {1}{\sqrt m}   \hspace{1em}  
        \\
        &= \frac{1}{m+1} \sum_{i=1}^{m}  \frac{2x}{M_{l_i-1} }  \, + \, \overset{s}{\sim} \frac {1}{\sqrt m}  
\end{align*}
\begingroup
\renewcommand{\thefootnote}{\fnsymbol{footnote}}
\footnotetext[1]{Among the subintervals containing $y$, some may have cutting points to the right of $y$, outside $(0,y)$.
}
\endgroup  
When $x=x_c$, the left-hand side of the above equals 1. 
Thus, we obtain the following lemma after simple algebra.
\begin{lemma} For any $y \in (0,1)$, the empirical c.d.f. of spacings within the subset $(0,y)$ is
\[
        \frac{1}{m+1} \sum_{i=1}^{m} \xi_{l_i}(x) = \frac{x}{x_c} 
        \, + \, \overset{s}{\sim} \frac {1}{\sqrt m}  
\]
for $x \leq x_c$. 
\end{lemma}
\noindent Hence, the number of partition points in $(0, y)$ is $2y/x_c$. Incorporating the fact that the number of partition points in $(0,1)$ is $2/x_c$, this implies 
\begin{theorem}  In the model of Kakutani interval splitting, the empirical c.d.f. of partition points satisfies      
\begin{equation}
                          F_n (x) = x \, + \, \overset{s}{\sim} \frac{1}{\sqrt{n}} \hspace{1em}  
\end{equation}
\end{theorem}

%\begin{align*}
%    a &= b + c 
%\begingroup\renewcommand{\thefootnote}{\fnsymbol{footnote}}\footnotemark[3]\endgroup \\
%    d &= e + f
%\begingroup\renewcommand{\thefootnote}{\fnsymbol{footnote}}\footnotemark[1]\endgroup \\    
%\end{align*}
%\begingroup
%\renewcommand{\thefootnote}{\fnsymbol{footnote}}
%\footnotetext[3]{First footnote explanation.}
%\endgroup 
 
%\begingroup
%\renewcommand{\thefootnote}{\fnsymbol{footnote}}
%\footnotetext[1]{First footnote explanation.}
%\endgroup  

%Furthermore, let
%\begin{align*}
%     &Y_{1}= X_{1}^2 - E X_{1}^2 \\
%     &Y_{2}= X_{2}^2 - E X_{2}^2 \\
%     &Y_{3}= X_{3}^2 - E_{1} X_{3}^2 \\
%     &\vdots \displaybreak[0]\\
%     &Y_{n}= X_{n}^2 - E_{n-2} X_{n}^2 \\
%\end{align*}
%It is easy to check that
%Furthermore, $X_{i}^2 = \xi_{i}^2(x)+E^2(\xi_{i}(x)) - 2\xi_{i}(x)E(\xi_{i}(x))$
\textbf{Remarks.} 
\begin{itemize}
  \item It is not necessary to assume the split point is uniformly distributed in Lemma~\ref{lem:log n}. This means that for a stationary process, the empirical partition points always approach the uniform distribution in the limit, irrespective of the distribution of the split points. 
  
  \item If and only if the mean spacing becomes stationary, in the sense that it converges to the same value in every subinterval $(0, y)$ for all $y$, the empirical c.d.f. of partition points converge to $x/x_c$ with the same fluctuation error. Dean and Majumdar (2002) and Janson and Neininger (2008) \cite{DeanMajumdar2002, JansonNeininger2008} gives sufficient conditions for the mean spacing to be stationary (inferable from their results). Therefore, their result can be extended to the equidistribution of empirical partition points. 
     
  %\item Further, we know that as long as mean spacing converges, in the sense that mean spacing converges to the same in any subinterval of size greater than the natural fluctuation $ \sim\! 1/\sqrt{n}$, the empirical c.d.f. of partition points is $x$. \cite{Janson and Neininger} have proved conditions for the mean spacing to be convergent (inferable from their statement (i) in Theorem 1.3). Therefore, the claim of theirs can be extended to c.d.f. of empirical partition points. 
     
  \item From the proofs of the theorems, the empirical c.d.f. of the spacings depends only on the distribution of the split points and is independent of the splitting scheme; the value of $\sum$ is invariant under permutations of the summands. This confirms Maillard and Paquette's (2016) conjecture \cite{MaillardPaquette2016}. 
       
  \item Our theorems provide the rate of convergence to the limiting distributions of spacings and partitions points. This may be an answer to another open question raised by Maillard and Paquette in the same paper. 
  
\end{itemize}

%Let $M_{1}^{(y)}, M_{2}^{(y)}, ..., M_{n_{y}}^{(y)}$ be all the intervals, in $(0, y)$, split during the process and %$\xi_{1}^{(y)}(x), \xi_{2}^{(y)}(x), ..., \xi_{n_{y}}^{(y)}(x)$ be the numbers the subintervals $\leq x$ generated by %those $n_{y}$ splittings.

%\[
%    \sum_{i=1}^{n^{(y)}} \xi_{i}(x) =  \sum_{i=1}^{n^{(y)}} E(\xi_{i}(x)) + \sim\! \sqrt n
%\]

\section{Distribution function of limiting spacings in general case} 
In this section, $F(x)$ related to $V_i$ is no longer assumed to be linear. We assume $F(x)$ is uniformly continuous, $F'(x)$ exists and uniformly positive. By~(\ref{E:empirical}), we have 
\begin{equation}\label{E: F cont}
    \frac{1}{N} \sum_{i=1}^{n} \xi_{i}(x) = \frac{1}{N} \sum_{i=1}^{n} F(x/M_{i-1}) 
    \, + \, \overset{s}{\sim} \frac {1}{\sqrt n}
\end{equation} where $N = n(b-1) + 1$. 
%where $\displaystyle N^{-1} \sum_{i=1}^{n} F(x/M_{i-1}) $ is continuous function of $x$ with error term of $\sim\! 1/\sqrt{N}$. 
For example, if $F(x)$ = $  a_1 x + a_2 x^2 $ then, the right-side of (\ref{E: F cont}) is 
\[ 
        \frac{a_1x}{1+n(b-1)} \sum_{i=1}^{n}   \frac 1 {M_{i-1}}       +  
        \frac{a_2x^2}{1+n(b-1)} \sum_{i=1}^{n} \frac 1 {M_{i-1}^2}
                                                             \, + \, \overset{s}{\sim} \frac {1}{\sqrt n}
\]
Note that although $\frac{1}{N} \sum_{i=1}^{n} F(x/M_{i-1})  $ is a continuous function of $x$, it does not necessarily converge as $n \to \infty $. If the convergence occurs  we call the fragmentation process is stationary. 

A key advantage of w.h.p. bound approach, compared with that of limit theorem, is that the union of many events holds also w.h.p. given each holds w.h.p. Specifically, if each of poly$(n)$ positions fails with exponentially small probability $\le \exp(-n^\alpha)$, then, by a union bound, all hold simultaneously w.h.p. also, since 
\[      
    1 - n^A \cdot e^{-n^\alpha} = 1 - e^{-n^\alpha + A\log n} \approx 1 - e^{-n^\alpha} \hspace{1em} \mbox{(for large $n$)}
\]
This property justifies calculus-like operations as long as the incremental amount $\Delta $ does not fall below natural fluctuation scale $ O(1/\sqrt{n})$.  
 
During the process of splitting, we always take the current largest interval, denoted $x_c$, as the unit. Let $u$ denote empirical d.f. of spacings at time $t$, i.e. 
\[
    u(x, t) = \frac { \mbox {\# spacings of length} \leq x}{N}
\] for $x < x_c$\,. 
%Let $N_{s}$ be the smallest natural number for which the length of all subintervals $\leq s$, i.e
%\[
%   N_{s} = min\{n \in Z^+: M_{n} \leq s \} 
%\]
%In other words, $N_{s}$ is the number of splittings immediately after all the subintervals are smaller than $s$ in the process. Let $t$ denote the time of this point.
By time $t + \Delta t$, all subintervals of length (normalized) $ > 1-\Delta s$ disappear (i.e. are split into shorter ones). %After rescaling, the interval $\left[0, 1-\Delta s \, \right ]$ is stretched to [0,1], mapping $1 - \Delta s$ to 1, and previous $x \cdot (1-\Delta s)$ is mapped to $x$. 
The total number of subintervals changes from $N$ to $N + \Delta N$. The proportion of subintervals with (previously) length $ > 1 - \Delta s$ is 
$$
                1 - u(1 - \Delta s, t) \approx u_x(1,t)\Delta s
$$
Thus, the number of subintervals split is 
\[
            N \cdot u_{x}(1,t)\Delta s
\]
Each split increases the number of subintervals by $(b - 1)$. This gives 
\begin{equation}\label{Delta N}
        \Delta N = (b - 1) \cdot  N \cdot u_x(1,t) \Delta s 
\end{equation}
The count of intervals $ \leq x\cdot (1-\Delta s)$ at time $t$, i.e.
                           $$ u(x\cdot (1-\Delta s), t) \cdot N $$
becomes the count of intervals $ \leq x$ $\big($after rescaling, the interval $\left[0, 1-\Delta s \, \right ]$ is stretched to [0,1], mapping $1 - \Delta s$ to 1, and previous $x \cdot (1-\Delta s)$ is mapped to $x\big)$.  Splittings during $\Delta t$ generate 
\[
        \frac { F(x) \Delta N} {b-1}
\]
additional intervals of length $<x$. Thus we have
\begin{align*}
    u(x, t + \Delta t) &= \frac{u(x(1-\Delta s), t)\cdot N + F(x) \Delta N /(b-1)}
                                 { N + \Delta N } \\[10pt]
    &\approx u(x(1-\Delta s), t) [1 - (b-1)u_{x}(1,t)\Delta s] + F(x)u_{x}(1,t)\Delta s 
                                                                    \hspace{2em}  \mbox{(by (\ref{Delta N}))}
                          \\[10pt]
                          &= u(x(1-\Delta s), t) - (b-1)u(x(1-\Delta s), t) u_{x}(1,t)\Delta s  + F(x)u_{x}(1,t)\Delta s
\end{align*}
i.e. 
\begin{align*}
  u(x, t + &\Delta t)  - u(x(1-\Delta s), t) \\
  &= - (b-1)u(x(1-\Delta s), t) u_{x}(1,t)\Delta s  + F(x)u_{x}(1,t)\Delta s
\end{align*}
With $\Delta s$ as measurement of $\Delta t$, some computation yields a first-order transport equation with a novel boundary feedback $u_x(1,t)$,
\/********************************************************************** 
\footnote{
One type of equation for fragmentation / renewal with boundary flux feedback
\[
     xu_x + u_t = S(t)(-\rho u + K(x)) 
\]
When $S(t) = u_x(1,t)$,  $\rho = (b-1)$ and $K(x) = F(x)$, it becomes our equation (\ref{E:partial}), i.e.
\[
     xu_x + u_t = u_x(1,t)(-(b-1)u + F(x)) 
\]
}
************************************************************************/ 
%\begin{align*}
%   x\frac {\partial u} {\partial x } + \frac {\partial u}{\partial t} &= - u(x, t)u_{x}(1,t) + u_{x}(1,t)L(x)/(b-1)
%\end{align*}
%i.e.
\begin{align}\label{E:partial} 
   x\frac {\partial u} {\partial x } + \frac {\partial u}{\partial t} 
   = u_{x}(1,t)\big(-(b-1)u(x, t)  +  F(x) \big)
\end{align}  
%The rate equations have not the first term, $x\cdot {\partial u} /{\partial x }$, of the above.%
%If $u$ converges the same distribution as time $t$ progresses, then ${\partial u(x,t)}/{\partial t} = 0$ and~(\ref{E:partial}) becomes
%\begin{equation}\label{ux}
%     {x} \frac {du}{dx} =  u_{x}(1,t) \big ( -(b-1)u(x) + F(x) \big)
%\end{equation} 
%This is a necessary condition for the stationary distribution of spacing. Kakutani's interval splitting process, where $u(x)= x$, satisfies above equation.

%Solution of~(\ref{ux}) is
%\[
%          u(x) = \frac {1}{x^a} \left( 1-a\int_0^1 \! t^{a-1}H(t) \, dt + a\int_0^x \! t^{a-1}H(t) \, dt\right)
%\] where $a = u_{x}(1)$ and $H(x) = L(x)/(b-1)$.

\textbf{Mean spacing stationarity  $\Rightarrow $ spacing distribution stationary.} Suppose the mean spacing (the first term below) is stationary (or equidistribution of partition points). We will show that $u$ is stationary (or ``stable'', as it is called in the area of PDEs), and provide its solution. 

Integration on both sides of~(\ref{E:partial}) yields
\begin{align}\label{E:int42}
   \int_0^1 \! x \frac {\partial u(x, t)} {\partial x } \, dx
    &  + \frac {\partial}{\partial t} \int_0^1 \! u(x, t)   \, dx   
    \\\nonumber
    & + u_{x}(1,t)\left[(b-1)\int_0^1 \! u(x, t)dx - \int_0^1 \! F(x)  \, dx \right]=0
\end{align}
By integration by parts, the mean spacing is 
\[
       \int_0^1 \! x \frac {\partial u(x, t)} {\partial x } \, dx  
              =  x \cdot u(x,t)  \Bigr|_0^1 - \int_0^1 \! u(x, t)   \, dx = \text {constant}
\]
$\iff$
\[   
       u(1,t) - \int_0^1 \! u(x, t)   \, dx = \text {constant} \hspace{2em}   
\] 
Since $u(1,t) = 1$,    
\[   
      \iff \int_0^1 \! u(x, t)   \, dx = \text {constant}  
\] 
and thus,  $$ \displaystyle \frac {\partial}{\partial t} \int_0^1 \! u(x, t)   \, dx = 0$$
It follows from (\ref{E:int42}) that $u_{x}(1,t)$ is a constant (let $u_{x}(1,t)=C$) and (\ref{E:partial}) becomes a first-order linear partial differential equation
\begin{align}\label{E:partial1}
   x\frac {\partial u} {\partial x } + \frac {\partial u}{\partial t} = -C(b-1)\cdot u + C\cdot F(x)
\end{align}
%The characteristic equations are
%\begin{align}\label{E:characteristic}
%              &\frac{dx}{dv}= x   % & & \text{($C=1$ because $\frac{dx}{dt}=1$ at $x=1$)}
%                                                                                              \notag \\
%              &\frac {dt} {dv} = 1  % & & \text{  (i.e. $dv=dt$  )}
%                                                                                           \\
%              &\frac {du}{dv} = -C\cdot u + C\cdot F(x)                           \notag
%\end{align}
%The first two equations of the above gives $x = C_{u} e^{t}$, where $C_{u}$ is a constant.
The characteristic equations of (\ref{E:partial1}) in the nonparametric form are
\begin{equation}\label{E:nonparam}
       \frac{dx}{x} = \frac{dt}{1} = \frac{du}{-C(b-1)\cdot u + C\cdot F(x)}
\end{equation}
It follows that $x=C_1 e^t$ where $C_1$ is an arbitrary constant, and
\begin{equation}\label{E:du/dx}
    \frac {du}{dx} + \frac{C(b-1)}{x}u = C\frac{F(x)}{x}
\end{equation}
\/********************************************************   
The solution of $du/dx + P(x)u=Q(x)$ is given by
\[
    u  = \Big(\int Q(x)e^{\int P(x)dx}dx + C_2 \Big) / e^{\int P(x)dx}
\] where $C_2$ and arbitrary constant.
**********************************************************/
The solution of (\ref{E:du/dx}) is 
\[
    u = \frac{1}{x^{C(b-1)}} \left (  C  \int_0^x F(s)s^{C(b-1)-1}ds +  a(xe^{-t})   \right)
\]
where $a(.)$ is an arbitrary function, $F(x)$ is defined for $x > 0$. Applying boundary condition $u(1,t)=1$, we find $a(\cdot)$ is a constant which must be 0 because $u$ is bounded at $x=0$ and $x^{-C(b-1)}$ diverge as $x \downarrow 0$. The solution of (\ref{E:du/dx}) turns out to be
\begin {equation}\label{uB}
    u = C x^{-C(b-1)}\int_0^x  F(s)s^{C(b-1)-1}ds
\end{equation} 
\/********************************************************   
% Try to find u'(1)
\begin {equation} 
    u = \frac{u_x(1)}{x^{u_x(1)(b-1)}}\int_0^x  F(s)s^{u_x(1)(b-1)-1}ds
\end{equation}

\begin {equation} 
    u'(x) = \frac{C \cdot F(x)}{x} - \frac{C(b-1)}{x}u(x)
\end{equation}
$\Rightarrow$
\[  u_x(1) = {C \cdot F(1)} - {C(b-1)} \hspace{1em} \mbox{(This is tautology and does not give $u_x(1)$, i.e. $C$)} \]
********************************************************************************/
where $C\, \big ( := u_x(1) \big) $ is defined above.  Then by $u(1,t)=1$, we have
\begin{equation}\label{stationary}
              1 = C \int_0^1  F(s)s^{C(b-1)-1}ds    %\mbox{\ \ \ (from $u(1)=1$)} 
\end{equation}
This gives the spacing c.d.f of a stationary sequential fragmentation process; for instance, in the case of ternary (uniform) splitting points, 
\[ 
        u(x) = \frac 1 2 x(3-x)
\] 
Thus, the following theorem is justified.
\begin{theorem}\label{limit mean spacing iff}
      In a sequential fragmentation, the spacing distribution is stationary if and only if the mean spacing is stationary. 
\end{theorem}

\begin{corollary}\label {E:div_space} Empirical distribution of partition points is uniform if and only if spacing's distribution is stationary. 
\end{corollary}

\textbf{Remarks.} 
\begin{itemize}

  \item  By Theorem \ref{limit mean spacing iff}, the results on the stationarity---sufficient conditions for stationarity---of the mean spacing obtained by Dean \& Majumdar (2002) and Janson \& Neininger (2008) \cite{DeanMajumdar2002, JansonNeininger2008} apply directly to that of the spacing distribution. 
  
  \item The equation \eqref{stationary} admits the unique solution $C = 1/(b-1)$ (see Appendix~\ref{q = 1}), independent of the distribution of the division points.
\end{itemize}

\section{Subtree size profile and fragmentation stationarity} 
A sequential interval splitting process can be encoded as a rooted tree, and vice versa. Each subinterval corresponds to a child subtree, with its length determining the subtree's size---for example, number of leaves \cite{DeanMajumdar2002, JansonNeininger2008, Sibuya&Itoh1987} or, in continuous random trees, total mass \cite{Aldous1991}. Tree profiles can be characterized by various properties, such as their height or their subtree size distribution, both can be expressed as interval fragmentation processes. While tree height has attracted much attention, and many interesting results have been obtained (see, for example, \cite{FlajoletSedgewick2009, Devroye1986})\,\footnote{From the perspective proposed by this paper, the height of a tree can be viewed as the steps for $y$-interval to reach the cutoff, $O(\log N)$, as determined by the shrinkage rate of the $y$-interval (see Appendix~\ref{k=log N}).}, our focus here is the \textit{subtree size profile}. 

For subtree size distribution, $m$-ary search trees \cite{Chauvin2011, ChauvinPouyanne2004, Mahmoud1992} and random recursive trees (Catalan trees) \cite{ChangFuchs2010, FlajoletSedgewick2009} are well studied. We give a general method to determine the subtree size profile for a broad class of random trees.

From (\ref{E:xi}), $f(x)\,dx $ = $F'(x)\,dx$ is the expected number of spacings with lengths in $(x, x+dx)$ from a split of $(0,1)$. For instance, in the Kakutani model $f(x) = 2$. In the case of 3 uniformly distributed split points, $f(x) = 6(1-x)$.

Let $g(x, y)\,dx$ denote the number of spacings with length in $(x, x+dx)$ at the moment when the current maximal spacing is evolves down to $y$.\,\footnote{After sufficiently steps $n$ of splitting, the actual number of this concentrates around the expected value, with deviations on the scale of natural fluctuation of order $1/\sqrt{n}$, with high probability.} An analysis similar to that in the previous section yields the evolution equation 
\begin{equation}
    g(x, y-dy)  = g(x,y) +  f\left(\frac {x}{y} \right) \frac {1}{y}\; g(y,y)\, dy \hspace{1em} \mbox{$(x < y)$} 
\end{equation}  
This leads to  
\begin{equation}
\nonumber      g(x,y) = g(x,1) + \int_y^1 f\left( \frac{x}{t} \right) \frac 1 t \; g(t,t)\,dt  \hspace{.5em} 
\begingroup
\renewcommand{\thefootnote}{$\ddagger$}
\footnote{In car parking problem, for car size $c$,  $y > c + x $, 
                \[   
                    g(x,y)  = \frac{1}{1-c} f\left(\frac{x}{1-c}\right) 
                           + \int_{y}^1 f\left(\frac {x}{t-c}\right) \frac {1}{t-c} g(t,t) dt 
                \]             
} 
     \addtocounter{footnote}{-1}     
\endgroup 
\end{equation}             
In the interval splitting process, the initial condition is given by $g(x,1) = f(x)$; this is the density of spacing immediately after the first split.  

Setting $x=y$ in the preceding equation yields an integral equation for the density of the largest spacing, defined as $H(y) := g(y,y)$, 
\begin{equation} \label{E: H(x)}
      H(y) = f(y) + \int_y^1 f\left( \frac{y}{t} \right) \frac 1 t H(t)dt  \hspace{.5em} 
\end{equation} 
\/********************************************************************************************************************

                                                    1    | 1                           1
  g(x,y) = 2 + 2C?_y^1  1/t^3 dt  =  2 + 2C (-1/2 -----) |    = 2 + 2C (-1/2 + (1/2) ----- 
                                                   t^2   | y                          y^2
 
                        C
         = 2 - C  +  -------
                       y^2
                       
  1/?N ~ ?y
                       
**********************************************************************************************************************/
In Kakutani model, $H(y) = Cy^{-2}$ where constant $C$ is determined by the total mass, the length of the unit interval $(0, 1)$. Thus, in $(0,y)$, $g(x,y) = H(y)$ (we already know $u'(x)$ is a constant function).  From the fact that the total mass is 1, 
\[
         1 = \int_0^y x \cdot H(y)\, dx %\ \Rightarrow \  H(y) \frac{1}{2}y^2 = 1 
         \ \Rightarrow \  C = 2 \ \Rightarrow \  H(y) = \frac {2} {y^2}
\]

\/******************************************************** 
In Kakutani model, $f = 2$. 
\begin{equation}\label{E: H(x)}
      H(x) = 2\int_x^1   \frac 1 t H(t)dt  \hspace{.5em} 
\end{equation} 
%                            /|  
%                          /  |         g(x,y) = C   F(x) = C·x
%                        /    |
%                       ·-----?    
%                          y
%
From Theorem \ref{E:th_spacing} we know $g(x,y)$ is a constant function of $x$. Let it be $C$ (distribution is $C\cdot x$). From 
\[
           \int_0^y C \cdot x \,dx = 1 \hspace{1em} \mbox{i.e. total mass, where $g(x,y)=C \cdot x$}
\] 

$\Rightarrow$ \  $\displaystyle  \frac C 2 y^2 = 1 \Rightarrow \ C = \frac 2{y^2} $  
********************************************************************************/ 
%In discrete context, splitting of $(0, l)$ produces $\displaystyle f({k}/{l})\cdot \frac {1}{l}$ spacings of size between $k$ and $k+1$ ($\displaystyle {1}/{l}$ corresponds to $dx$). % Therefore, after an iteration, $g(k)$ gains $ f(\frac {k}{l}) (\frac {1}{l}) h_l(l)$. 
In the context of trees, let $n$ denote the size of the tree (number of leaves). The total mass is $n$. Encoding a sequential interval splitting process as a tree structure, the subtrees of size $k$ corresponding precisely to the largest remaining subintervals of lenght $k$. Conversely, in a tree, labelling all roots of size-$k$ subtrees with the same time stamp produces a time-indexed family of subtrees that can be viewed as a sequential interval fragmentation process. Then, the number of subtrees with size $k$, i.e. the number of largest subintervals, is given by $n^{-1} H(k/n)$ and others' are $n^{-1} g(i/n, k/n) $. For instance, in the case of binary search tree, 
\[
         \frac{1}{n} H(k/n)= \frac {2n} {k^2}
\]

In a tree encoding of the splitting process, the number of subtrees of size $k$ corresponds to the number of largest subintervals of normalized length $k/n$. This is given by $n^{-1} H(k/n)$; all other subtrees are counted by $n^{-1} g(i/n, k/n)$. For instance, in a binary search tree,
\[
    \frac{1}{n} H(k/n) = \frac{2n}{k^2}.
\]

\/******************************************************** 
To find $C$, note the total mass is $n$
\begin{align*}
        n=   \sum_{i=1}^k  i\cdot g_n(k)   &=  n^2 \sum_{i=1}^k  \frac i n \cdot g_n(k) \frac 1 n  
        \\ 
                                           &= n^2 \int_0^{k/n}  x\cdot \frac {2} {(k/n)^2}  dx
        \\                                           
                                           &=n^2 \frac {2} {(k/n)^2} \cdot \frac 1 2 {\left( \frac{k}{n} \right)}^2 
\end{align*}
********************************************************************************/ 
\/********************************************************
In discrete context, splitting of $(0, l)$ produces $\displaystyle f(\frac {k}{l})\cdot (\frac {1}{l})$ spacings of size between $k$ and $k+1$ ($\displaystyle {1}/{l}$ corresponds to $dx$). 
% Therefore, after an iteration, $g(k)$ gains $ f(\frac {k}{l}) (\frac {1}{l}) h_l(l)$.

In car parking problem, splitting of splitting of $(0, l)$ produces $ f(\frac {k}{l-1}) (\frac {1}{l-1})$ spacings of size between $k$ and $k+1$. After an iteration of parking on spaces of length $l$, $g(k)$ gains $ f(\frac {k}{l-1}) (\frac {1}{l-1}) h_l(l)$.

Let splitting process undergo from max interval of $n$ down to $k$, starting from the interval $[0,n]$. Let
\[
     \ h_n(k) = \delta (k-n)
\] 
where $\delta$ is the Dirac delta function. 
At each splitting step, the recurrence $h$, the count of spacings of size $k$,  is given by 
\begin{align*} 
    & h_{n-1}(k) = h_n(k) + f(\frac {k}{n})  \frac {1}{n}  h_{n}(n)
                           \mbox{, ~ splitting $(0, n)$}  \\
    & h_{n-2}(k) = h_{n-1}(k) + f(\frac {k}{n-1})  \frac {1}{n-1}  h_{n-1}(n-1)
                                 \mbox{, ~ splitting $(0, n-1)$} \\
    & ... \\
    & h_{l-1}(k) = h_{l}(k) + f(\frac {k}{l})  \frac {1}{l}  h_{l}(l)
                                 \mbox{, ~ splitting $(0, l)$} \\
    & h_{l}(k) = h_{l+1}(k) + f(\frac {k}{l+1})  \frac {1}{l+1}  h_{l+1}(l+1)
                                 \mbox{, ~ splitting $(0, l+1)$} \\
    & ...  \\
    & h_{k}(k) = h_{k+1}(k) + f(\frac {k}{k+1})  \frac {1}{k+1}  h_{k+1}(k+1)
                                 \mbox{, ~ splitting $(0, k+1)$}
\end{align*}
As a result, for $k \le l$,
\[
    h_{l}(k) = h_n(k) + \sum_{j=n}^{l+1} f(\frac {k}{j}) (\frac {1}{j}) h_{j}(j)
\]
$h_n(k)$ is negligible if $k \ll n$, i.e.
\[
    h_{l}(k) = 0 + \sum_{j=n}^{l+1} f(\frac {k}{j}) (\frac {1}{j}) h_{j}(j)
\]
This is about expectation. We need to consider the real time (random) number $\widetilde{h}(k)$, i.e. 
\[
   \widetilde{h}_{l}(k) = 0 + \sum_{j=n}^{l+1} \widetilde{f}\left(\frac {k}{j} \right) 
                                                \frac {1}{j}\cdot 
                                               \widetilde{h}_{j}(j)
         \mbox {~~~ ~$ \left( \displaystyle \sum_{k=0}^n \sum_{j=n}^{l+1} \widetilde{f}(\frac {k}{j}) \hdots  (k/n)^{s-1}
                       \right)  
                     $ 
               }
\]
$\widetilde{f}(\frac {k}{j}) \to f(\frac {k}{j})$ so long as $\frac {1}{j} \gg \frac{1}{\sqrt{n}}$
It is observed that $\widetilde{h} \gg 1$ for $l \ll \sqrt{n}$. Around $l=\sqrt{n}$, $\widetilde{h}_{l}(k) \approx $ continuous function but $|h|$ is very small.

The summation can be written as an integral,  
\[
             h_y(x)  = h_n(x) + \int_{y}^1 f(\frac {x}{t}) \frac {1}{t} h_t(t) dt \hspace{1em} 
\]  
where $x=\frac {k}{n}$, $t=\frac {j}{n}$, $\frac{1}{j} = \frac{1}{n} \frac{n}{j} = \frac{1}{n} \frac{1}{t}=dt/t$, $h_1(t) := h_j(j)$ and $h(x,y) := h_l(k)$. In particular, for $l=k$, $h_1(x)$ satisfies the following integral equation,
\[
    h_1(x-dt)  = \delta(1-x)  + \int_{x}^1 f(\frac {x}{t}) \frac {1}{t} h_1(t) dt
                      \mbox {~~ ($1/n$ omitted)}
\]
\[\Big(
     h_1(x)  = \delta(1-x) + \int_{y}^1 f(\frac {x}{t}) \frac {1}{t} h_1(t) dt
\Big) \]

\[\Big(
     h(x,y)  = \delta(1-x) + \int_{y}^1 f(\frac {x}{t-c}) \frac {1}{t-c} h_1(t) dt
     \mbox { ~~ for car parking of car size $c$, $y > c + x$
           }
\Big) \]

When $n$ is extremely large compared to $k$, $h_n(k)$ is very small compared to $h_1(x)$. In this case, $h_n(x)$ is negligible. That is to say, for $x$ is small ($k$ is extremely small compared to $n$)

\[
   \widetilde{h}_1(x-T)  =  \int_{x}^1 f(\frac {x}{t}) \frac {1}{t} h_1(t) dt
\]

\[
    h_1(x)  =  \int_{x}^1 f(\frac {x}{t}) \frac {1}{t} h_1(t) dt
\]

\[\Big(
    h_1(x)  =  \int_{x+c}^{1} f(\frac {x}{t-c}) \frac {1}{t-c} h_1(t) dt
\Big)\]
********************************************************************************/ 
For $m$-ary search tree with splitting probability density $f(x)= m(m-1)(1-x)^{m-2}$, (\ref{E: H(x)}) becomes 
\begin{equation}\label{eqn: y(x)}
       H(y) = f(y) +  m(m-1) \int_{y}^1 (1 - \frac {y}{t} )^{m-2} \frac {1}{t} H(t) dt
\end{equation}
In this Volterra integral equation, if the non-homogeneous term $f(x)$ remains bounded as $x \to 0$, it may be neglected when analysing the singular behavior of $y(x)$ for small $x$. 
Differentiating $m-1$ times with respect to $y$ eliminates the integral (see Appendix~\ref{I(m-1)}) and gives 
\begin{equation}\label{eq: Cauchy-Euler} 
   H^{(m-1)}(y)  =  (-1)^{m-1} m! \cdot    \frac {1}{y^{m-1}} H(y)
\end{equation} 
a Cauchy-Euler differential equation. Chern et al.\ (2002) \cite{Chern2002} first found the connection from $m$-ary search tree to Cauchy-Euler differential equation. General solution of (\ref{eq: Cauchy-Euler}) is
\[
  H(y) =   C_1 y ^{\lambda_1} + C_2 y ^{\lambda_2} + \hdots + C_{m-1} y ^{\lambda_{m-1}} 
                                                                               \, + \, \overset{s}{\sim}  \frac {1}{\sqrt{n}}   
            \hspace{1em}   
   \mbox{($n$: times of splitting)}  \hspace{-4.8em}
\]
i.e.
\[
  H(y) =  C_1 \left(\frac{1}{y} \right)^{2} +\; C_2 \left(\frac{1}{y} \right) ^{z_2 } + \hdots 
                                                                               +\; C_{m-1} \left(\frac{1}{y} \right) ^{z_{m-1} }
           \, + \, \overset{s}{\sim}  \frac {1}{\sqrt{n}}
\]
$\lambda_i$ is the $i^{\text{th}}$ solution of
\/******************************************************
\[
      \frac{1}{(-1)^{m-1} m(m-1)(m-2) ...2 } \cdot \lambda(\lambda - 1)...
      (\lambda - m+2 )  = 1
\]
i.e.
*******************************************************/
\[
        \lambda(\lambda - 1)... (\lambda - m+2 )  = {(-1)^{m-1} m! }
\]
Letting $z = -\lambda $, the above equation is transformed equivalently to
\begin{equation}\label{E: psi(z)}
         z(z+1)... (z + m-2) =  m!
\end{equation}  
$z_1 (=2), z_2, . . ., z_{m -1}$ denote the $m - 1$ roots of (\ref{E: psi(z)}) in nonincreasing order of real parts and roots with positive imaginary parts listed before their conjugates. 

\textbf{Example: }  For $m=3$,  
\[
  H(y) =  C_1 \left(\frac{1}{y} \right)^2 +\;  C_2 \cdot y^{3}  
\]
From the initial condition $H(y) = 6(1-y)$ which gives $H(1)=0$ and $H'(1)=-6$, we obtain 
\[
           H(y) =  \frac{6}{5}   \left(\frac{1}{y^2}  -  y^{3}  \right)
\]
In the context of a ternary search tree, the number of subtrees of size $k$ is given by  
\[
           \frac{1}{n} H(k/n)   \approx  \frac{6}{5}  \frac{n}{k^2}  
\]
neglecting the term of $ k^3/n^4 $ for $k \ll n$. To the best of our knowledge, the subtree-size profile is known only for random binary search tree \cite{Dennert2010} \cite{Fuchs2008} and random recursive (binary) tree \cite{Feng2008} \cite{Devroye1991}.

The relationship between the roots of (\ref{E: psi(z)}) and the ``phase change'' in $m$-ary search trees has been studied extensively (e.g. \cite{Mahmoud1992, Hwang2003, ChauvinPouyanne2004}). Here, we interpret the phase transition as a measure of stationarity -- mentioned in previous sections -- in the associated interval fragmentation process. In other words, the stationarity of spacings and of the distribution of partition points in a fragmentation process can be studied through the analysing the roots of a Cauchy-Euler differential equation, or integral equation (\ref{eqn: y(x)}) in general. The results obtained in the existing literature on the ``phase change'' apply to the stationarity of the $m$-ary fragmentation process; for example, when $m \le 26$ the $m$-ary fragmentation is stationary. %\,\footnote{Rigorously speaking, this may only be a necessary condition for $u(x,t)$ to be stationary. It is unclear how sufficient it is. Moreover, formal analysis of the relationship between $H(y)$ and $u(x,t)$ is required. As this paper aims to introduce novel methodological frameworks based on the new conceptual perspective, we do not go into the technical detail here and leave these questions open.} 
 The method indicated here is general and applies to arbitrary splitting probability density $f(x)$.

\section{Heavy-tailed distribution of split points} In the previous section $f(x)$ denotes a continuous density on $(0,1)$. We now consider the discrete regime, so that the initial unit interval is $[1, n]$, where the split point distribution is given by a heavy-tailed probability mass function $(p_i)_{1\le i \le n}$. Specifically, without loss of generality, assume
\[
         p_i \propto i^{-\alpha}, \qquad 1 < \alpha < 2,
\] 
An important example is the fringe-tree distribution in Catalan trees,\,\footnote{
$ \displaystyle p_i = \frac{1} {2\sqrt{\pi} i^{3/2}} \big( 1 + \frac{1}{8(i-1)}  + \cdots  \big ) $
}  where $p_i \propto i^{-3/2}$. The discrete regime requires different analytical tools from the continuous case. In the current literature, subtree size profiles of Catalan trees are typically studied using generating function method, yielding a few expectation-based results % rather than a complete characterization of their behaviour 
\cite{FlajoletSedgewick2009,ChangFuchs2010}.
 
%Let the split process start with the unit interval $(0,n)$ with total mass $n$. 
As the fragmentation process is also called \textit{mass fragmentation}, we may refer to an interval as a ``mass''. In one step, a mass of size $n$ is likely to be cut into a very small piece of size $i$ and a much larger of $n-i$, with probability $p_i$. Strictly speaking $p_i$ depends the parent's length. For instance, if $p_i \asymp  1/i^\alpha$, then, for a parent of length $n$, the exact values are  
$$
    \frac{1}{\zeta_n(\alpha)1^\alpha},\;\;    \frac{1}{\zeta_n(\alpha)2^\alpha},\;\; 
                \cdots \;\;  \frac{1}{\zeta_n(\alpha)i^\alpha},\;\;  \cdots
                                                                            % , \;  \frac{1}{\zeta_n(\alpha)(n/2)^\alpha} 
$$
where $$\zeta_n(\alpha) := \zeta (\alpha) - O(1/n^{\alpha - 1})$$ 
However, except for small mass of order $O(1)$, the difference between $\zeta_n(\alpha)$ and the constant $\zeta(\alpha)$ may be ignored, so that $p_i$ can be treated as nearly constant independent of $n$. For example, for $ \alpha = {3/2}$, the difference between $\zeta_n(3/2)$ and $\zeta(3/2)$ is as small as $O(1/\sqrt{n})$. Henceforth, we will use $\zeta(3/2)$ or even $\zeta$, for $\zeta_n(3/2)$ to simplify the notation. 

For convenience, we label the smaller cut pieces ``black'' masses and the larger as ``white'' masses. So in the greedy splitting scheme, each step only a white mass is cut; black masses accumulate and remain untouched until white masses shrink to comparable size. 

\textbf{Recurrence for $H_t$\,.} \,Note that the subintervals all originate from the historically largest masses. Let $H_t$ be the number of largest masses at time $t$, with $H_0 = 1$. After $t$ steps, all masses of size $> n-t$ have been cut, and 
\begin{equation}\label{eq:Ht_rec}
    H_t = p_1 H_{t-1} + p_2 H_{t-2} + \cdots + p_{t-1} H_1 + p_t H_0,
\end{equation}
is a linear difference equation of unbounded order. It is practically difficult to solve a general linear difference equation of unbounded order.\,\footnote{For discussion of the solutions for general linear difference equations with variable coefficients, see, e.g., \cite{Mallik1998} and references therein. }\,  Fortunately (\ref{eq:Ht_rec}) is a restricted type in which the coefficients form a probability density, and we can solve it. Iterating the recursion (\ref{eq:Ht_rec}) yields 
\[
\begin{pmatrix} H_t \\ H_{t-1} \\ \vdots \\ H_0 \end{pmatrix}
= C^{\,t}(p) \begin{pmatrix} H_0 \\ 0 \\ \vdots \\ 0 \end{pmatrix},
\]
where $C(p)$ is the companion matrix built from $(p_1, p_2, \dots)$. For instance, at time $1$ 
\[
\begin{pmatrix}
    H_1     \\  H_0     \\  0       \\ \vdots  \\  0  
\end{pmatrix}
 =
\begin{pmatrix}
    p_1 & p_2 & p_3 &\cdots   & p_{n/2} \\
    1 & 0 & 0 & \cdots & 0 \\
    0 & 1 & 0 & \cdots & 0 \\
    \vdots & \vdots & \ddots & \ddots & \vdots \\
    0 & 0 & \cdots & \ \ 1 & 0    
\end{pmatrix} 
\begin{pmatrix}
    H_0     \\   0     \\   0       \\   \vdots  \\   0  
\end{pmatrix}     
= C(p) 
\begin{pmatrix}
    H_0     \\   0     \\   0       \\   \vdots  \\   0  
\end{pmatrix}     
\]
at time $2$ 
\[
\begin{pmatrix}
    H_2     \\   H_1   \\  H_0  \\   \vdots  \\   0  
\end{pmatrix}
 =
\begin{pmatrix}
    p_1 & p_2 & p_3 &\cdots   & p_{n/2} \\
    1 & 0 & 0 & \cdots & 0 \\
    0 & 1 & 0 & \cdots & 0 \\
    \vdots & \vdots & \ddots & \ddots & \vdots \\
    0 & 0 & \cdots & \ \ 1 & 0    
\end{pmatrix} 
\begin{pmatrix}
    H_1     \\ H_0     \\   0       \\ \vdots  \\ 0  
\end{pmatrix}     
= C^{\,2}(p) 
\begin{pmatrix}
    H_0     \\   0     \\   0       \\   \vdots  \\   0  
\end{pmatrix}     
\] 
%\/******************************************************
By iteration,  
\[
\begin{pmatrix}
    H_t  \hspace{0.8em}   \\  H_{t-1} \\   H_{t-2} \\   \vdots \hspace{0.7em}  \\   H_0 \hspace{0.7em} 
\end{pmatrix}
 = C^{\,t}(p)
\begin{pmatrix}
    H_0     \\   0     \\   0       \\   \vdots  \\   0   
\end{pmatrix}     
\]
%*******************************************************/
From the Perron-Frobenius theorem \cite{Seneta1981}, we know 
\[
                  \lambda_1 =1 > | \lambda_i | \quad\mbox{(for $i \neq 1$ )}
\]
%$C(p)$ can be written as $ V J V^{-1} $
%\[
%      J^t \rightarrow  \operatorname{diag}(1, 0, 0, \cdots, 0)
%\]
Thus, after a large number \big(but negligible relative to $n$, e.g., $ O(\log n)$\big) steps, $C^{\,t}(p)$ approaches the rank-one trivial matrix associated with $\lambda_1 = 1$, namely\,\footnote
{The eigenvalues are distinct generically and $C(p)$ is diagonalizable: 
$ C(p) = V\operatorname{diag}(1,\, \lambda_2^t,\, \lambda_3^t,\, \cdots \; )V^{-1}$.
}
\[
           C^t(p) \rightarrow  V\operatorname{diag}(1,\, \lambda_2^t,\, \lambda_3^t,\, \cdots \; ) V^{-1}
                  \rightarrow   \operatorname{diag}(e_{11},\, 0,\, 0,\, \cdots \; )  
\]
Hence, $H_t$ stabilizes until a later critical time.

\medskip
\noindent
\textbf{Example: $\alpha = 3/2$}.
For $p_i \propto i^{-3/2}$, one finds (see Appendix~\ref{H_t})
\[
    H_t = \frac{\zeta(3/2)}{4\sqrt{M}},
\]
where $M$ is the distance from the origin.  
One may picture the fragmentation as rapidly breaking a single thick interval into many---on the order of $\sqrt{n}$\,--- thin intervals of similar size, each of amount $O(1/\sqrt{n})$. The largest subintervals number $\displaystyle \zeta/ (4\sqrt{n})$. In parallel, fringe intervals are created with the distribution $p_i \propto 1/i^{3/2} $ for intervals $(1, i)$ and remain untouched until the critical point.

\textbf{Critical points.} While the number of white largest intervals stay essentially constant, the number of black intervals are increasing by 
$$
  \displaystyle   \frac {O(1)} { i^{3/2}} 
$$
for the intervals $(1, i)$. A critical point occurs when $H_t$ begins to increases suddenly, because those historically black intervals start to dominate the regime of largest interval while the white masses vanish. Technically, at this point the largest black masses and largest white masses are equal in size and quantity, i.e.,   
\[\displaystyle
        \frac{C_M}{M^{3/2}}= H_t = \frac {\zeta  } {  4\sqrt{M}    }  \quad \Rightarrow \quad   C_M  = \frac {\zeta  }  4 M
\]   
By mass conservation,  
\[ n = C_M\sum_{i=1}^M i\cdot p_i = C_M\sum_{i=1}^M i \cdot \frac 1 {i^{3/2}}
\]  
hence $ 
        C_M 2\sqrt{M} \approx n
      $. 
Solving the two equations gives
\[ 
   M = M_c  \approx 0.84 n^{2/3}\; , \quad \displaystyle C_{M_c} = \zeta  M_c /4 \approx 0.55n^{2/3} \hspace{.5em}
\begingroup
\renewcommand{\thefootnote}{$\ddagger$}
\footnote{$ \displaystyle M_c = \left( \frac{2}{\zeta} \right)^{2/3} n^{2/3}$\,, \hspace{.0em} 
          $ \displaystyle C_{M_c} = \frac{\zeta}{4} \left( \frac{2}{\zeta} \right)^{2/3} n^{2/3}$
} 
     \addtocounter{footnote}{-1}     
\endgroup                  
\]  
and the distribution of spacings
\begin{equation}
            p_i = \frac {C_{M_c}}{i^{3/2}}   =  O \left( \frac{n^{2/3}}{i^{3/2}} \right)   \quad  1 \le i \le  M_c        
\end{equation}  
%Chang and Fuchs~\cite{ChangFuchs2010},  in the context of Catalan tree, identified a critical point where the mean of $X_{n,k}$ (normalized number of subtrees with size $k$) equals 1; their threshold has the same order $n^{3/2}$, and the exact constant is immaterial. For $k \sim n$ they obtained $X

We next find $H_k$ for $M_c' \le k \le M_c$, where $k$ is the distance from the origin 1, and $M_c'$ is the next critical point. Consider subintervals of length $i$, $k \le i \le M_c$. Each single subintervals of these generates the same $H_t = \zeta /(4\sqrt{k})$,\,\footnote{$H_t$ of an interval is independent of its length; the length only determines the position of its critical points.}   so 
\[
         H_k  = \frac{\zeta}{4\sqrt{k}} \;  \sum_{i=k}^{M_c} \frac{C_{M_c}}{i^{3/2}} \approx 
         \frac{\zeta}{2{k}}C_{M_c} \left( 1 - \sqrt{\frac k {M_c} } \right) \;
%         \quad \quad 0.84 {M_c}^{2/3} \le k \le M_c  
          \quad \quad M_c' \le k \le M_c  
\]
In the black regime, increments always have distribution of the form 
\[
         \frac {C_k} {i^{3/2}}
\]
with the amplitude $C_k$ increasing over time. The next critical point $M_c'$ occurs when the number of largest black equals $H_{M_c'}$\,:
\[
    \frac{C_{M_c'}}{M_c'^{3/2}} = H_{M_c'} =  \frac{\zeta}{2{M_c'}}C_{M_c} \left( 1 - \sqrt{\frac {M_c'} {M_c} } \right)  
\]
Mass conservation gives
\[
       \sum_{i=1}^{M_c'} i \cdot \frac {C_{M_c'}} {i^{3/2}}  = n
\]
\/*********************************************************
The first equation (A) gives

    C_Mc' = ?/2(Mc') C_Mc · { 1 -  ? (Mc'/Mc ) } Mc'³?²        
      => 
    C_Mc' = ?/2 C_Mc ?Mc'  
                                 
The second (B) gives
           Mc'                       ___         
    C_Mc'  ?    1/?i  =>    C_Mc' 2 ?Mc'  = n   
          i=1                                  
           
 Mc' = 1/? · n/C_Mc ~ n¹?³ 
 
     C_Mc' = ??/2  ?n ?C_Mc  = ??/2 * ?n  ?(? Mc/4) = ?/4 ?(n M_c) ~ n???
                       ~~~~
                         \
                        C_Mc = ? Mc/4 
                   1      
     Mc' = 4n/?² -----  ~ n¹?³  
                  Mc                           
***********************************************************/
Solving yields \[
\displaystyle  M_c' =  \frac{4n} {\zeta^2 M_c}  = O( n^{1/3} )\; , \quad \displaystyle C_{M_c'} = \frac{\zeta}{4}\sqrt{n} \sqrt{M_c}  =  O( n^{5/6}  ) \]  
Thus 
\begin{equation} 
                p_i = \frac {C_{M_c'}}{i^{3/2}} = O \left( \frac{n^{5/6}}{i^{3/2}} \right) \quad  1 \le i \le  M_c' 
\end{equation}
From this spacing distribution one can determine successive critical points, until $k$ is sufficiently small that the normalized spacing distribution is $O(1)$.

\textbf{Reverse recurrence for small $k$.}
In the small-$k$ regime, we use a backward recurrence, starting from $k=1$, as illustrated below. 

\noindent $\mathbf{k=1}$: $n_1(1) = H_1 = n$

\noindent$\mathbf{ k=2 }$: 
\vspace{-1em}
\begin{align*}
                 n_1(1) &= n_2(1) + p_1\cdot n_2(2) 
               \\
                 n \hspace{1.em}  &= n_2(1) + \ 2\cdot n_2(2)   \quad \mbox{(mass conservation)}
\end{align*} 
Equivalently, 
\[
   A_2 \begin{pmatrix}
                 n_2(1)  \\[4pt]
                 n_2(2)  
      \end{pmatrix} =  
    \begin{pmatrix}
      n_1(1) \\[4pt]
      n
    \end{pmatrix} ,  \quad   
        A_2 = 
        \begin{pmatrix}
            1 & p_1 
          \\[4pt]
            1 & 2
        \end{pmatrix}              
\]
 
\noindent $\mathbf{k=3}$: 
\vspace{-1em}
\begin{align*}
                  n_2(1) &= n_3(1) +  p_1\cdot n_3(2) +  p_2\cdot n_3(3) 
               \\
                  n_2(2) &= \hspace{5.2em}    n_3(2)  +  p_1\cdot n_3(3) 
               \\               
                     n \hspace{1.em}  &= n_3(1) + \ \  2\cdot n_3(2) + 3 \cdot n_3(3)    
\end{align*}  
i.e.
\[
   A_3    \begin{pmatrix}
                 n_3(1)  \\[4pt]
                 n_3(2)  \\[4pt]
                 n_3(3)
      \end{pmatrix} =  
    \begin{pmatrix}
      n_2(1) \\[4pt]    
      n_2(2) \\[4pt]
      n
    \end{pmatrix}
     ,    \quad   
      A_3 = 
        \begin{pmatrix}
           1   &  p_1   & p_2   \\[4pt]
           0   &  1     & p_1   \\[4pt]        
           1   &  2     & 3
        \end{pmatrix}         
\]    
and $p_1 = p_2$. 

Generally case: for $k \ge 2$, 
\begin{equation}\label{E: A_k}
   A_k  
     \begin{pmatrix}
                 n_k(1)    \\[4pt]
                 n_k(2)    \\[4pt] 
                 \vdots    \\[4pt]
                 n_k(k-1)  \\[4pt]
                 n_k(k)
      \end{pmatrix} =     
      \begin{pmatrix}
                 n_{k-1}(1)    \\[4pt]
                 n_{k-1}(2)    \\[4pt] 
                 \vdots        \\[4pt]
                 n_{k-1}(k-1)  \\[4pt]
                 n
      \end{pmatrix} 
                                              : =  
                                                \begin{pmatrix}
                                                     \boldsymbol{b}  \\[4pt]          
                                                            n                      
                                               \end{pmatrix}    
\end{equation}
\[               
  A_k = 
        \begin{pmatrix}
            1     &   p_1     &   p_2     &   p_3      &  \cdots   &      p_{k-1}   \\
            0     &   1       &   p_1     &   p_2      &  \cdots   &      p_{k-2}   \\
            0     &   0       &   \ddots  &    \ddots  &  \ddots   &       \vdots   \\
          \vdots  &  \vdots   &   \ddots  &   \ddots   &   p_1     &      p_2       \\
            0     &   0       &   \cdots  &   0        &  1        &      p_1       \\
            1     &   2       &  3        &   \cdots   &  k-1      &      k                
        \end{pmatrix} 
\]
where $p_i$ satisfy\,\footnote{In the case of Catalan tree, $ 
                \displaystyle p_i = \frac{C_{i-1} C_{k-i} } {C_k}, \quad C_n = \frac{1}{n+1} \binom{2n}{n} $ 
}
\[
   p_i = p_{\,k-i}, \quad 1 \le i \le k-1
\] 

Unwinding the recurrence gives
\[
     \begin{pmatrix}
                 n_k(1)    \\[4pt]
                 n_k(2)    \\[4pt] 
                 \vdots    \\[4pt] 
                 n_k(k)
      \end{pmatrix}    
     = \prod_{j=1}^{k-1}      
                        \begin{pmatrix}
                            A_{j}^{-1} &     0 
                                                     \\[4pt]
                            0          &   I_{k-j}  
                        \end{pmatrix}      
                                       \begin{pmatrix}
                                                   n     \\[4pt]
                                                   n    \\[4pt] 
                                                   \vdots    \\[4pt] 
                                                   n 
                                        \end{pmatrix}      
\]
with $A_1^{-1} = I_1$. 

We leave open finding the asymptotic $H_k := n_k(k)$; it is conjectured 
$H_k =  O\left( \displaystyle \frac{n}{k^\alpha} \right) $. 
\/**************************************************************************************************************************
\textbf{Estimate $H_k$}.  Let $ \boldsymbol{u}=\left(p_{k-1}, p_{k-2},...,p_1\right)^{\mathrm{T}}$, $\boldsymbol{v}  = \left(1,2, ..., k-1\right)^\mathrm{T}$. Write
\[  A_k = 
          \begin{pmatrix}
               \boldsymbol{U}       &  \boldsymbol{u} 
                              \\[4pt]
               \boldsymbol{v}^{\mathrm{T}}     &  k
          \end{pmatrix}  ,   \quad  U \in \mathbb{R}^{(k-1)\times(k-1)}
          %\,  \mbox{upper-triangle with 1's on the diagonal}
\]       
where  $\boldsymbol{U}$ is an upper-triangular matrix with 1's on the diagonal and $p_1, p_2, ..., p_{k-2}$ on the first, second, ..., $k-2)$-th super-diagonals, respectively. For any right-hand side (see (\ref{E: A_k}))
\[
                            \begin{pmatrix}
                                 \boldsymbol{b}  \\[4pt]          
                                        n                                                                                     
                             \end{pmatrix}     \quad \mbox{(here  $\boldsymbol{b} = \boldsymbol{n}_{k-1}$)}   
\]
Schur complement gives 
\[        
          n_k(k) = S^{-1} \left( n - \boldsymbol{v}^{\mathrm{T}} \boldsymbol{U}^{-1} \boldsymbol{b} \right), \quad
          S := k - \boldsymbol{v}^{\mathrm{T}} \boldsymbol{U}^{-1} \boldsymbol{u}
\]
\[
      \boldsymbol{U} ^{-1} = \left(\boldsymbol{I} + \boldsymbol{N} \right)^{-1}  
      = \boldsymbol{I} - \boldsymbol{N} + O(\boldsymbol{N}^2) % - \boldsymbol{N}^3 + \cdots
\] 
Thus 
\begin{equation}
    \boldsymbol{v}^{\mathrm{T}} \boldsymbol{U}^{-1} \boldsymbol{b} = 
                                   \boldsymbol{v}^{\mathrm{T}} \boldsymbol{b} 
                                -  \boldsymbol{v}^{\mathrm{T}}  \boldsymbol{N} \boldsymbol{b}  
                                +  \boldsymbol{v}^{\mathrm{T}}  O(\boldsymbol{N}^2) \boldsymbol{b} 
%                                -  \boldsymbol{v}^{\mathrm{T}}  \boldsymbol{N}^3 \boldsymbol{b} 
%                               + \cdots
\end{equation}
where
\[               
  \boldsymbol{N} =  
        \begin{pmatrix}
            0     & \  \ p_1   &  \ p_2       &  \quad  \cdots   &          \ p_{k-3}  &  \ p_{k-2}    \\
            0     &   0       &   \ p_1       &  \quad  \cdots   &          \ p_{k-4}  &  \ p_{k-3}    \\
            0     &   0       &   \ddots      &   \ddots         &                     &  \vdots       \\
          \vdots  &  \vdots   &   \ddots      &   \ddots         &   \quad  \ p_1      &  \ p_2        \\
                  &           &   \cdots      &   \quad          &   \quad   0         &  \ p_1        \\
            0     &   0       &   \cdots      &   \quad          &                     &    0            
        \end{pmatrix}        
        = 
        \begin{pmatrix}
            0     & \  \ p_1   &  \ p_2       &  \quad  \cdots   &     \quad p_{3}   &  \quad  p_{2}    \\
            0     &   0       &   \ p_1       &  \quad  \cdots   &     \quad p_{4}   &  \quad  p_{3}    \\
            0     &   0       &   \ddots      &   \ddots         &                   &  \quad  \vdots   \\
          \vdots  &  \vdots   &   \ddots      &   \ddots         &    \quad  \ p_1   &  \quad p_2        \\
                  &           &   \cdots      &   \quad          &    \quad   0      &  \quad p_1        \\
            0     &   0       &   \cdots      &   \quad          &                   &  \quad  0            
        \end{pmatrix}  
\]

\[      \boldsymbol{I} - \boldsymbol{N} =  
        \begin{pmatrix}
            1     & \  -p_1   &  -p_2         &   \cdots   &   \quad -p_3      &  -p_{2}    \\
            0     &   1       &   -p_1        &   \cdots   &   \quad -p_{4}    &  -p_3      \\
            0     &   0       &   \ddots      &   \ddots   &                   &  \vdots    \\
          \vdots  &  \vdots   &   \ddots      &   \ddots   &   \quad  -p_1     & -p_2       \\
            0     &   0       &   \cdots      &   \quad    &   \quad   1       & -p_1       \\
            0     &   0       &   \cdots      &   \quad    &                   &  1            
        \end{pmatrix} 
\]
\[
    \boldsymbol{v}^{\mathrm{T}} \boldsymbol{N} \boldsymbol{b}
    = 
        \begin{pmatrix}       1     &  \     2       &  \  3        &   \cdots   &  k-1      \end{pmatrix} 
        \begin{pmatrix}
            0     & \  \ p_1   &  \ p_2       &  \quad  \cdots   &     \quad p_{3}   &  \quad  p_{2}    \\
            0     &   0       &   \ p_1       &  \quad  \cdots   &     \quad p_{4}   &  \quad  p_{3}    \\
            0     &   0       &   \ddots      &   \ddots         &                   &  \quad  \vdots   \\
          \vdots  &  \vdots   &   \ddots      &   \ddots         &    \quad  \ p_1   &  \quad p_2        \\
                  &           &   \cdots      &   \quad          &    \quad   0      &  \quad p_1        \\
            0     &   0       &   \cdots      &   \quad          &                   &  \quad  0            
        \end{pmatrix}       
        \begin{pmatrix}
              n_{k-1}(1)          \\
              n_{k-1}(2)          \\
              \vdots              \\
                                  \\
              n_{k-1}(k-2)        \\
              n_{k-1}(k-1)    
        \end{pmatrix}         
\]    

\[
        \begin{pmatrix}       1     &  \     2       &  \  3        &   \cdots   &  k-1      \end{pmatrix} 
      \left(
             \begin{array}{l}
                p_1 + p_2 + \cdots + p_k / 2 + \cdots + p_3 + p_2 \\
                p_1 + p_2 + \cdots \qquad \cdots +   p_4 +    p_3 \\
                \vdots \\
                p_1 + p_2 \\
                p_1 \\
                0
            \end{array}
      \right)
      : = sum
\]
Let     
\[
         v_j \approx 2 \sum_{i=1}^{\lfloor (k-j)/2 \rfloor} p_i
\]

\[            sum =  \sum_{j=1}^{k-1} j\cdot v_j \approx 2\sum_{j = 1}^{k-1} j\sum_{i=1}^{\lfloor (k-j)/2 \rfloor} p_i
\]
 
\[
        \sum_{i=1}^m  \frac{1}{i^\alpha} \sim \zeta(\alpha)  - \int_m^{\infty}   \frac{1}{x^\alpha} \,dx
                 = \zeta(\alpha) - \frac{m^{1-\alpha}}{\alpha - 1}
\] 
\[ 
          \sum_{i=1}^{\lfloor (k-j)/2 \rfloor} \frac{1}{i^\alpha} 
                       \sim \zeta(\alpha) - \frac{  (k-j)^{1-\alpha}  }   {2^{1-\alpha}(\alpha - 1)}
\]
This gives 
\[
       sum \approx 2 \sum_{j=1}^{k-1} j \left[ \zeta (\alpha) - C(k-j)^{1-\alpha} \right]
             = 2\zeta (\alpha) \sum_{j=1}^{k-1} j   - 2C  \sum_{j=1}^{k-1}   j(k-j)^{1-\alpha}  
\]

\[
         \sum_{j=1}^{k-1}   j(k-j)^{1-\alpha}  \sim k^{2-\alpha} \cdot \int_0^1 x(1-x)^{1-\alpha}\,dx \sim \Theta(k^{2-\alpha})
\]
Thus, 
\[
           sum \sim  \zeta(\alpha) k^2 - Ck^{2-\alpha} 
\]
The leading term is  $\zeta(\alpha) k^2 $. 
 
By mass conservation, the zeroth term is 
\[ 
      \boldsymbol{v}^{\mathrm{T}} \boldsymbol{b}  =  \sum_{j=1}^{k-1} j\cdot n_{k-1}(j) = n
\]
Therefore
\[
     \left( n - \boldsymbol{v}^{\mathrm{T}} \boldsymbol{U}^{-1} \boldsymbol{b} \right) 
              \approx \boldsymbol{v}^{\mathrm{T}}  \boldsymbol{N} \boldsymbol{b}
              \sim \zeta(\alpha) 2n
\]

For the first correction $\boldsymbol{v}^{\mathrm{T}}  \boldsymbol{N} \boldsymbol{b}  $, 
\[
     (\boldsymbol{Nb})_i = \sum_{j = 1}^{k-1-i} p_j n_{k-1}(i+j) 
\]
$        \Rightarrow   $
\[
     \boldsymbol{v}^{\mathrm{T}}  \boldsymbol{N} \boldsymbol{b} = \sum_{i=1}^{k-1} i \sum_{j = 1}^{k-1-i} p_j n_{k-1}(i+j)
                       = \sum_{j = 1}^{k-2} p_j \sum_{i=1}^{k-j-1} i \cdot n_{k-1}(i+j)   
\]
where $i+j \le k-1$. Therefore,  
\[
     \left( n - \boldsymbol{v}^{\mathrm{T}} \boldsymbol{U}^{-1} \boldsymbol{b} \right) 
              \approx \boldsymbol{v}^{\mathrm{T}}  \boldsymbol{N} \boldsymbol{b}  
\]      
%  v^{T}N b = n/k^2 { 
%                  p_1 (1 + 2  +  3 + ...                    + (k-2) )
%              +   P_2 (1 + 2  +  3 + ...               + (k-3)      )
%              +   P_3 (1 + 2  +  3 + ...           + (k-4)          )
%                  .
%                  .
%                  .
%                p_{k-3)(1 + 2                                       )
%             +  p_{k-2)(1                                           )      
%         }
%    = { 1*sum_1^{k-2} p_i + 2*sum_1^{k-3} p_i + 3*sum_1^{k-4} p_i +  (k-2)*p_1  }
% 
%             k-2
%    =  (1/2) ?    p_i (k - 1 - i)(k - i)  n/k^2
%             i=1
% 
%     ~ O(k^2)
%  
% 
% 
% 
%
\[
          S := k - \boldsymbol{v}^{\mathrm{T}} \boldsymbol{U}^{-1} \boldsymbol{u} 
             = k - \boldsymbol{v}^{\mathrm{T}} (\boldsymbol{I} - \boldsymbol{N})  \boldsymbol{u} 
     \approx   k - \sum_{j=1}^{k-1}  j \cdot p_j + \boldsymbol{v}^{\mathrm{T}} \boldsymbol{N} \boldsymbol{u}
\]
 
\[
         T_k = \sum_{j=1}^{k-1}  j \cdot p_j \sim \frac{c}{2-\alpha}k^{2-\alpha} , \quad (2 - \alpha < 1)
\] 

\[
    \boldsymbol{v}^{\mathrm{T}} \boldsymbol{N} \boldsymbol{u}
    = 
        \begin{pmatrix}       1     &  \     2       &  \  3        &   \cdots   &  k-1      \end{pmatrix} 
        \begin{pmatrix}
            0     & \  \ p_1   &  \ p_2       &  \quad  \cdots   &     \quad p_{3}   &  \quad  p_{2}    \\
            0     &   0       &   \ p_1       &  \quad  \cdots   &     \quad p_{4}   &  \quad  p_{3}    \\
            0     &   0       &   \ddots      &   \ddots         &                   &  \quad  \vdots   \\
          \vdots  &  \vdots   &   \ddots      &   \ddots         &    \quad  \ p_1   &  \quad p_2        \\
                  &           &   \cdots      &   \quad          &    \quad   0      &  \quad p_1        \\
            0     &   0       &   \cdots      &   \quad          &                   &  \quad  0            
        \end{pmatrix}       
        \begin{pmatrix}
              p_{1}           \\
              p_{2}           \\
              \cdot              \\
              p_{k/2}       \\
              \cdot               \\     
              p_2                \\
              p_1   
        \end{pmatrix}         \sim   \Theta(k^{3-2\alpha})
\]    
 
So the only place the history $A_1, ..., A_{k-1}$ appears is via  $\boldsymbol{b} = \boldsymbol{n}_{k-1}$ in the scalar $\boldsymbol{v}^{\mathrm{T}} \boldsymbol{U}^{-1} \boldsymbol{b} $
%
%              n
% So n_k(k) ~ ----    This method ignore n_{k-1}(i) ; that might be why it fails to obtain n/k?
%              k
%
***********************************************************************************************************************/
\clearpage
\appendix 

\section{\centering} \label{Azuma-Heoffding inequality}
%\/******************************************* commented out: ...  
\noindent \textbf {Concentration inequality (Azuma-Heoffding inequality)}
\begin{proof} 
Let $ \displaystyle S_m = \sum_{i=1}^m X_{l_i} $ \,  $l_1 < l_2 < \hdots  <  l_m$. 
\begin{align*}  \vspace{1em}              
        \mathbb{E} \left[  e^{t S_m } \right] 
   &=  \mathbb{E} \left[  e^{t ( X_{l_1} + X_{l_2} + \hdots + X_{l_m} ) } \right]
                   \\[10pt]
   &= \mathbb{E} \left[ \mathbb{E} \left(  e^{t ( X_{l_1} + X_{l_2} + \hdots + X_{l_m} ) }  
                                               \,\middle|\,  M_{l_m - 1}, ..., M_1, M_0 \right) \; 
                \right] 
                   \\[10pt]
   &= \mathbb{E} \Big[\, e^{t S_{m-1}} \mathbb{E} \left(  e^{t X_{l_m}  }  
                                               \,\middle|\,  M_{l_m - 1}, ..., M_1, M_0 \right)\; 
                 \Big]  
\end{align*}
Applying $ \mathbb{E} [ X_{l_m}  \mid  M_{l_m - 1}, ..., M_1, M_0 ] = 0$, and $|X_i| < c$,  by Heoffding's Lemma,
\begin{align*} 
           \mathbb{E} \Big[  e^{t X_{l_m}  }  \mid  M_{l_m - 1}, ..., M_1, M_0 \Big] \le e^{t^2c^2/2}   
\end{align*}
Thus, \vspace{-1em}
\begin{align*} 
     \mathbb{E} \left[  e^{t S_m } \right] \le \mathbb{E} \left[ e^{t S_{m-1}} \right] e^{t^2c^2/2}       
\end{align*}
Iteration obtains 
\vspace{-1em}
\begin{align*} 
     \mathbb{E} \left[  e^{t S_m } \right] \le e^{mt^2 c^2/2}       
\end{align*}
For $\lambda > 0$, by Markov's inequality\,\footnote
{$\mathbb{P}(X \ge a) \le \displaystyle \frac{\mathbb{E}[X]}{a}$. Thus 
                  $\mathbb{P}(S_m \ge \lambda) = \mathbb{P}(tS_m \ge t\lambda) = \mathbb{P}(e^{tS_m} \ge e^{t\lambda})   
                         \le  \displaystyle e^{-t\lambda}{\mathbb{E} \big( e^{tS_m} \big) }   $
}
\vspace{-0.5em}  
\begin{align*} 
      \mathbb{P}(S_m \ge \lambda) \le e^{-t \lambda} \mathbb{E} \left[  e^{t S_m } \right] \le       
      \exp \Big (-t \lambda + m\frac{t^2c^2}{2}    \Big)
\end{align*}
Optimize $t = \lambda/(mc^2)$, 
\vspace{-1.5em}
\begin{align*} 
      \mathbb{P}(S_m \ge \lambda) \le  \exp \Big (- \frac{\lambda^2}{2mc^2}  \Big)
\end{align*} 
Combine with $\mathbb{P}(S_m < -\lambda)$, 
\begin{align*} 
      \mathbb{P}(|S_m| \ge \lambda) \le  2\exp \Big (- \frac{\lambda^2}{2mc^2}  \Big)
\end{align*}
\end{proof}
%*******************************************/
\/******************************************* commented out: ...  
\noindent Prove that   
      $S_m = X_{l_1}  + \hdots + X_{l_{m-1} } + X_{l_m}$, $l_1 < l_2 < \hdots < l_{m-1}  <  l_m$,  is a martingale with respect to $\{ M_{l_m - 1}, ..., M_1, M_0 \}$
\begin{proof} Given $ \mathbb{E} [ X_{l_m} \mid M_{l_m - 1}, ..., M_1, M_0 ] = 0$,  
\begin{align*}  
    \mathbb{E} \left[ S_m  \mid M_{l_m - 1}, ..., M_1, M_0 \right]     
   &=\mathbb{E} \left[ S_{m-1} + X_{l_m} \mid M_{l_m - 1}, ..., M_1, M_0 \right]  
    \\
   &= S_{m-1} + \mathbb{E} [ X_{l_m} \mid M_{l_m - 1}, ..., M_1, M_0 ] 
   \\
   &= S_{m-1} 
\end{align*}  
\vspace{-1em}  
\end{proof} 
*******************************************/

\clearpage
 
\section{\centering} \label{k=log N}
Consider the unit interval $(0,1)$, and fix a point $y \in (0,1)$. Perform a sequential splitting process as follows. At each step, retain the subinterval containing $y$, referred as $y$-interval, to further split and discard the other. Repeat this procedure until the size of the retained subinterval is smaller than $x_c$. 

Assume $y \ge 1/2$; otherwise take $(y, 1)$. We show $k = O( \log N )$  w.o.p., where $k$ is the times of splitting and $ N = 1/x_c \asymp n $. 
 
\begin{proof}
Let $I_{i-1}$ be the $ith$ retained subinterval encompassing the point $y$, $i=0, 1, 2, ...$, the total shrinkage rate $R_k$ after $k$ splittings is 
\[
               R_k =  r_k \cdot r_{k-1} \hdots r_1 
\] 
where $r_k$ is one step shrinkage rate on $I_{k-1}$ given $\{I_i \}_0^{k-1} $. This holds 
\begin{align*}  
    \mathbb{E} \left[ R_k   
     \right]
    &= \mathbb{E}  \Big[ \mathbb{E} \left[ r_k \cdot R_{k-1} \,\middle|\, I_{k-1}, ..., I_1, I_0  \right] \, \Big]  
     = \mathbb{E}  \Big[  R_{k-1}\cdot \mathbb{E} \left[ r_k  \,\middle|\, I_{k-1}, ..., I_1, I_0  \right] \, \Big]
    \\ 
    & \le  \mathbb{E} \Big[  R_{k-1}\cdot \frac34  \Big] =  \frac34 \mathbb{E} \Big[  R_{k-1} \Big]  
      \le  \left( \frac34 \right)^2  \mathbb{E} \Big[  R_{k-2} \Big] \le  \hdots  \le  \left(\frac34 \right)^k    
\end{align*}  
where $ \displaystyle \frac34$ is the upper bound of the expected one-step shrinkage rate, achieved when $y$ is the middle point of the parent subinterval $I_{i-1}$. Now let $K$ represent the stopping time at which $y$-interval is first cut to a subinterval smaller than $x_c$, namely, 
\[
     R_{K}  \le \frac{x_c} y \le  R_{K-1} 
\]     
To estimate $K$, by Markov inequality, we have
\[
     \mathbb{P} \left(R_K > \frac{x_c} y \right) \le y\frac {\mathbb{E}[R_K]}{x_c} \le    N \left(\frac34 \right)^K 
\]
This implies $K$ =   $ O\left( \log N \right)$ w.o.p. In deed, for any $A > 0$, so long as  
\[ \displaystyle 
    K \ge  \frac{(1 + A)}{\log (4/3)}  \log N, 
\]    
\vspace{.3\baselineskip}
its failing probability $\mathbb{P}(R_k > x_c/y)$ is as small as $O(N^{-A})$.  
\end{proof}

\begin{proposition} \label{K/n = +- sg}
\[
          \frac{K} n = \pm_{sg} \frac{1}{\sqrt{n}}   
\]  
\end{proposition}

\begin{proof}  

Let $k_n$ be a number such that $ N \left(3/4 \right)^{k_n} = 1 $. For some $C_0$ and $c_0$, it holds
\begin{equation}\label{Ce^-cs}
         \mathbb{P} \left(K > k_n + s \right)  \le  \mathbb{P} \left(R_{k_n + s} > \frac{x_c} y \right) 
         \le \left( \frac 3 4 \right)^s = C_0 e^{-c_0 s }
\end{equation}
We aim to prove, for some $C$, $c$
\[
      \mathbb{P} \left(\frac K n > \frac{\lambda} {\sqrt{n}}  \right)  \le C e^{-c\lambda^2}
\]
To this end, we distinguish three cases for $\lambda$. 

  \textbf{Case 1.} $ \displaystyle \lambda \lesssim k_n /\sqrt n$. Since $k_n = O(\log n)$ and thus $ \displaystyle  \lambda = O\left( k_n / \sqrt n \right) \rightarrow 0$,  it holds that 
\[
     \mathbb{P} \left(\frac K n > \frac{\lambda} {\sqrt{n}}  \right) \le 1 \le C_1 e^{-c_1\lambda^2} 
\]
for some $C_1$ and $c_1$. 

\textbf{Case 2.} $ 2k_n /\sqrt n \, \le  \, \lambda  \, \le \, \sqrt n$. $ 2k_n /\sqrt n \, \le  \, \lambda$ implies 
\[                       
                       k_n   \le \frac{\lambda \sqrt n} 2  \;\; \iff \;\; 
                   \lambda \sqrt n  -  k_n   \ge \lambda \sqrt n - \frac{\lambda \sqrt n} 2    \;\; \iff \;\; 
                   \lambda \sqrt n  -  k_n   \ge  \frac{\lambda \sqrt n} 2
\]
Let $s = \lambda \sqrt n  -  k_n $.   $\eqref{Ce^-cs}$ yields 
\begin{align*}  
        \mathbb{P} \left(\frac K n > \frac{\lambda} {\sqrt{n}}  \right)  
        & \le C_0 e^{-c_0(\lambda \sqrt n  -  k_n) }  
        \\
        & 
         \le  
               C_0 e^{-c_0 (\lambda \sqrt n )/2 }   \qquad (\lambda \sqrt n  -  k_n   \ge  \frac{\lambda \sqrt n} 2)
       \\[4pt]              
        &\le  
               C_0 e^{-c_0 (\lambda^2)/2 }     \qquad (\lambda \le \sqrt n)
        \\
        &
          = C_2 e^{-c_2 (\lambda^2) }     \qquad (C_2:=C_0)
\end{align*}
Then, taking 
\[
       C = \max \{C_1, C_2\}, \quad c = \min \{c_1, c_2\}
\]
we have, for $\lambda \le \sqrt n$, 
\[
        \mathbb{P} \left(\frac K n > \frac{\lambda} {\sqrt{n}}  \right) \le C e^{-c\lambda^2}
\]
 
\textbf{Case 3.}  $\lambda > \sqrt n$.
In this case,  $\lambda \sqrt n > n$
Therefore, 
\[
    \mathbb{P} \left(\frac K n > \frac{\lambda} {\sqrt{n}}  \right)  = 
    \mathbb{P} \left(  K >  \lambda \sqrt n    \right) < \mathbb{P} \left(  K >  n  \right) = 0 \le C e^{-c\lambda^2}
\] 
Finally, we obtain, for any $\lambda$, 
\[         
         \mathbb{P} \left(\frac K n > \frac{\lambda} {\sqrt{n}}  \right) \le C e^{-c\lambda^2}
\]
Therefore 
\[
          \frac{K} n = \pm_{sg} \frac{1}{\sqrt{n}}   
\]
\end{proof}

\clearpage

\section{\centering} \label{I(m-1)}    
Let  $ \displaystyle         I(x) =  m(m-1) \int_{x}^1 (1 - \frac {x}{t} )^{m-2} \frac {1}{t} y(t) dt $
%\/******************************************************** 
\begin{align*}
       I'(x)  &= m(m-1)(-1) (1-x\frac{1}{x})^{m-2}y(x)\frac{1}{x}
       + m(m-1)(m-2) \int_{x}^1 (1 - \frac {x}{t} )^{m-3} (-1/t) \frac {1}{t} y(t) dt
              \\
              &= (-1)m(m-1)(m-2) \int_{x}^1  (1 -   \frac {x}{t} )^{m-3} \frac {1}{t^2} y(t) dt \\
       I''(x) &= (-1)^2m(m-1)(m-2)(m-3) \int_{x}^1  (1 - \frac {x}{t} )^{m-4}  \frac {1}{t^3} y(t) dt
             \\ 
              &\quad \vdots \\  
      I^{(m-2)}(x)   &= (-1)^{m-2} m(m-1)(m-2)(m-3)\hdots(m-(m-1)) \int_{x}^1  (1 - \frac {x}{t} )^{m-m}
                       \frac {1}{t^{m-1}} y(t) dt
                        \\
                     & = (-1)^{m-2} m! \int_{x}^1 \frac {1}{t^{m-1}} y(t) dt
                          \\
   I^{(m-1)}(x)  & =  (-1)^{m-1} m! \cdot    \frac {1}{x^{m-1}} y(x)
\end{align*}
%********************************************************************************/

\clearpage 
\section{\centering} \label{q = 1}   
Set, $q = C(b-1)$. Noticing that 
            \[
                  F(s) =  \mathbb{E} \left[ 
                                              \sum_{j=1}^b \mathbf{1}_{  \{ V_{j} < x \}  }  
                                     \right]
            \]              
we obtain 
\begin{align*}  
      \int_0^1 F((s) s^{q-1} ds & = \mathbb{E}\left[
                                                    \sum_{i=1}^b \int_{V_i}^1 s^{q-1} ds  
                                             \right]      
      \\
      & = \frac{1}{q} \mathbb{E}\left[ 
                                       \sum_{i=1}^b  \left(1 - V_i^q \right)
                               \right]      
      \\
      & =  \frac{1}{q}\left(b - \mathbb{E}  \sum_{i=1}^b  V_i^q  \right)
\end{align*}
Thus, the equation \eqref{stationary} 
\[ 
              1 = C \int_0^1  F(s)s^{C(b-1)-1}ds    %\mbox{\ \ \ (from $u(1)=1$)} 
\]
is equivalent to
\[
    1 =  \frac{1}{b-1}\left(b - \mathbb{E}  \sum_{i=1}^b  V_i^q  \right)
\]
i.e. 
\[
    \mathbb{E}  \sum_{i=1}^b  V_i^q  = 1
\]
This function of $q$ is strictly decreasing in $q$ (because $\{V_i\}_1^b$ is a genuine fragmentation and $0 \le V_i < 1$). Thus, $\mathbb{E}  \sum_{i=1}^b  V_i = 1 $ implying $q = 1$.  Consequently 
\[
        C = \frac{1}{b-1}
\]

\clearpage 
\section{\centering} \label{H_t}    
$ \displaystyle H_t = \frac{\zeta_n(3/2)}{4\sqrt{M}}$, for large $M$ which is the distance from the origin $1$.
\begin{proof}  
At time $t$, the number of interval $(1, n-t)$ is $H_t$. Let $H_t(k)$ denote the number of interval $(1, n-t-k)$, $k \ge 1$,  Then 
\begin{align*} 
        H_t(1) &= p_2H_{t-1} +  p_3H_{t-2} +  \cdots + p_{t}   H_{1} + p_{t+1} H_0 
        \\[10pt]
        H_t(2) &= p_3H_{t-1} +  p_4H_{t-2} +  \cdots + p_{t+1} H_{1} + p_{t+2} H_0     
        \\       
        \vdots
        \\
        H_t(i) &= p_{i+1}H_{t-1} +  p_{i+2}H_{t-2}   + \cdots + p_{t+i-1} H_{1} + p_{t+i} H_0 
        \\       
        \vdots                
\end{align*} where $H_0 = 1$. 
By the constancy of $H_{t}$, we have 
\begin{align*} 
        H_t(1) & \approx (p_2  +  p_3  +  p_4 \cdots ) H_t 
        \\[10pt]
        H_t(2) & \approx (p_3  +  p_4  +  p_5 \cdots ) H_t 
        \\       
        \vdots
        \\
        H_t(i) & \approx (p_{i+1}  +  p_{i+2}  +   p_{i+2}  \cdots)H_t
        \\  
        \vdots                
\end{align*}
and therefore
\[         
      H_t + \sum_{k=1}^M H_t(k) \approx H_t \sum_{k=0}^M \left(1 - \sum_{i=1}^k p_i \right) 
\]  
Let $p_i = \propto 1/i^\alpha$. Then, since $M$ is large, 
\[         
      H_t + \sum_{k=1}^M H_t(k) \approx   H_t\, \frac{1}{\zeta_n(\alpha)(\alpha - 1)(2-\alpha)} M^{2-\alpha}
\]  
When $\alpha = 3/2$ 
\[         
      H_t + \sum_{k=1}^M H_t(k) \approx   H_t \frac{4}{\zeta_n(3/2)}  \sqrt{M}
\]  
From mass conservation,  $H_t + \sum_{k=1}^M H_t(k) = 1 - o(1) $,\,\footnote{
          When $t  = O ( \log n )$, fraction $\displaystyle O\Big( \frac {\log n}{n} \Big)$ of intervals go black. } 
 we obtain, \[ \displaystyle H_t = \frac{\zeta_n(3/2)}{4\sqrt{M}} \] 

\/********************************************************************************  
\section{\centering} \label{appendix:e}    

Let $b_j = n_{k-1}(j) $ 
Mass conservation: 
\begin{equation}\label{n=}
            n =\sum_{j=1}^{k-1} j \cdot b_j = \sum_{j=1}^{m-1} j \cdot b_j + \sum_{j=m} ^{k-1} j \cdot b_j 
\end{equation} 
The second summation, 
\[
        \sum_{j=m} ^{k-1} j \cdot b_j  = \sum_{i=0}^{k-1-m} (i + m)\cdot b_{i+m} 
        = \sum_{i=1}^{k-1-m} i\cdot b_{i+m} + m  \sum_{i=0}^{k-1-m}  b_{i+m} 
\]
Plugging this into (\ref{n=}), we have 
\[
  n -  \sum_{j=1}^{m-1} j \cdot b_j = \sum_{i=1}^{k-1-m} i\cdot b_{i+m} +  m  \sum_{i=0}^{k-1-m} b_{i+m}                      
\]
We obtain, 
\[
  \sum_{i=1}^{k-1-m} i\cdot b_{i+m}  =  n -  \sum_{i=1}^{m-1} i \cdot b_i -   m \sum_{i=m}^{k-1} b_i 
                                     =  n -  m\left(\sum_{i=1}^{m-1} \frac{i}{m} \cdot b_i + \sum_{i=m}^{k-1} b_i \right)
\]
%              m\left(\sum_{i=1}^{m-1} \frac{i}{m} \cdot b_i 
Given the condition that $n_{k-1}(i)$ is not concentrating, 

********************************************************************************/

\end{proof}

\end{document}